\documentclass{article}
\usepackage{amsmath,amsthm,amssymb,amscd}
\theoremstyle{plain}

\newfont{\cyrrlarge}{wncyr10}
\def\Sha{\mbox{\cyrrlarge Sh}}

\newfont{\cyrr}{wncyr5}

\newtheorem{thm}{Theorem}[section]
\newtheorem{prop}[thm]{Proposition}
\newtheorem{lemma}[thm]{Lemma}
\newtheorem{cor}[thm]{Corollary}

\theoremstyle{definition}

\newtheorem{defn}[thm]{Definition}
\newtheorem{remark}[thm]{Remark}

\newcommand{\Z}{{\bf Z}}
\newcommand{\Q}{{\bf Q}}

\newcommand{\Zp}{\Z_p}
\newcommand{\bG}{\mathbb G}
\newcommand{\Zpx}{\Z_p^\times}
\newcommand{\Qp}{\Q_p}

\newcommand{\Fp}{{\bf F}_p}

\newcommand{\Qpbar}{\overline{\Q}_p}
\newcommand{\Qinf}{\Q_\infty}
\newcommand{\Qn}{\Q_n}
\newcommand{\Kinf}{K_\infty}

\newcommand{\inj}{\hookrightarrow}
\newcommand{\surj}{\twoheadrightarrow}
\newcommand{\maps}{\longrightarrow}
\newcommand{\lra}{\longrightarrow}

\newcommand{\da}{\downarrow}

\newcommand{\bH}{{\bf H}}
\newcommand{\ve}{\varepsilon}
\newcommand{\w}{\omega}
\newcommand{\wn}{\omega_n}
\newcommand{\tomega}{\tilde{\omega}}

\newcommand{\twnp}{\tomega_n^+}

\newcommand{\twnsubtwoe}{\tomega_{n-2}^\ve}
\newcommand{\twnm}{\tomega_n^-}
\newcommand{\twnjp}{\tomega_{n,j}^+}
\newcommand{\twnjm}{\tomega_{n,j}^-}
\newcommand{\twne}{\tomega_n^\ve}
\newcommand{\twnminuse}{\tomega_n^{-\ve}}
\newcommand{\twnje}{\tomega_{n,j}^\ve}

\newcommand{\twnmp}{\tomega_n^\mp}
\newcommand{\wnsubone}{\omega_{n-1}}
\newcommand{\wnp}{\omega_n^+}

\newcommand{\wnsubtwop}{\omega_{n-2}^+}

\newcommand{\wnsubtwoe}{\omega_{n-2}^\ve}
\newcommand{\wnm}{\omega_n^-}

\newcommand{\wne}{\omega_n^\ve}
\newcommand{\wme}{\omega_m^\ve}
\newcommand{\wnminuse}{\omega_n^{-\ve}}

\newcommand{\wnpm}{\omega_n^\pm}
\newcommand{\wnmp}{\omega_n^\mp}
\newcommand{\zn}{\zeta_n}

\newcommand{\lambdapm}{\lambda_E^\pm(\Kinf/K)}

\newcommand{\mupm}{\mu_E^\pm(\Kinf/K)}
\renewcommand{\H}{{\mathcal H}}
\renewcommand{\P}{{\mathcal P}}
\newcommand{\OO}{{\text{O}}}
\newcommand{\p}{{\mathfrak p}}
\newcommand{\cO}{{\mathcal O}}
\newcommand{\Xwnp}{\wnp}
\newcommand{\Xwnm}{\wnm}
\newcommand{\Xwne}{\wne}

\newcommand{\Xwnsubtwoe}{\wnsubtwoe}
\newcommand{\Xwnsubtwop}{\wnsubtwop}

\newcommand{\Xwnpm}{\wnpm}
\newcommand{\Xwnmp}{\wnmp}
\newcommand{\Xwnminuse}{\wn^{-\ve}}

\DeclareMathOperator{\ord}{ord}

\DeclareMathOperator{\Hom}{Hom}

\DeclareMathOperator{\Tam}{Tam}
\DeclareMathOperator{\Gal}{Gal}
\DeclareMathOperator{\Aut}{Aut}
\DeclareMathOperator{\Tr}{Tr}
\DeclareMathOperator{\Sel}{Sel}
\DeclareMathOperator{\coker}{coker}
\DeclareMathOperator{\im}{im}
\DeclareMathOperator{\cores}{cor}
\DeclareMathOperator{\rk}{rk}
\DeclareMathOperator{\cork}{corank}

\newcommand{\dual}{\wedge}
\newcommand{\Trn}{\Tr^n_{n-1}}
\newcommand{\kTn}{\ker(\Trn)}
\newcommand{\ckTn}{\coker(\Trn)}
\renewcommand{\char}{\text{char}}
\newcommand{\hE}{\widehat{E}}

\newcommand{\Sn}{S_n}

\newcommand{\Sne}{S_n^\ve}
\newcommand{\Sme}{S_m^\ve}
\newcommand{\Snme}{S_n^{-\ve}}
\newcommand{\Snsubone}{S_{n-1}}

\newcommand{\Sinfn}{S^{\Gamma_n}}
\newcommand{\Xn}{X_n}
\newcommand{\Xnsubone}{X_{n-1}}
\newcommand{\Xinfn}{X_{\Gamma_n}}
\newcommand{\Xinfm}{X_{\Gamma_m}}
\newcommand{\Xinfnsubone}{X_{\Gamma_{n-1}}}
\newcommand{\Xnpm}{X_n^\pm}
\newcommand{\SnT}{S_n(T)}
\newcommand{\SnTzero}{S_n^0(T)}
\newcommand{\SnTstrict}{S_{n,\Sigma}(T)}

\newcommand{\HinfjpmT}{{\bf H}^1_\pm(K_{\infty,\p_j},T)}
\newcommand{\HnjpmT}{H^1_\pm(K_{n,\p_j},T)}

\newcommand{\HnT}{H^1(K_{n,\p},T)}
\newcommand{\HneT}{H^1_\ve(K_{n,\p},T)}
\newcommand{\HnpT}{H^1_+(K_{n,\p},T)}
\newcommand{\HnmT}{H^1_-(K_{n,\p},T)}
\newcommand{\HinfpmT}{{\bf H}^1_\pm(K_{\infty,\p},T)}
\newcommand{\HnpmT}{H^1_\pm(K_{n,\p},T)}
\newcommand{\HnmpT}{H^1_\mp(K_{n,\p},T)}
\newcommand{\Ln}{\Lambda_n}
\newcommand{\Lnsubone}{\Lambda_{n-1}}
\newcommand{\Lnsubtwo}{\Lambda_{n-2}}
\newcommand{\Ld}{\Lambda^d}
\newcommand{\Lnd}{\Lambda_n^d}
\newcommand{\Lmd}{\Lambda_m^d}

\newcommand{\Lpmn}{\Lambda / \Xwnpm \Lambda}

\newcommand{\Len}{\Lambda / \Xwne \Lambda}
\newcommand{\Lminusen}{\Lambda / \Xwnminuse \Lambda}

\newcommand{\fgn}{\widehat{E}(K_{n,p})}
\newcommand{\fgm}{\widehat{E}(K_{m,p})}
\newcommand{\fgne}{\widehat{E}^\ve(K_{n,\p})}
\newcommand{\fgnp}{\widehat{E}^+(K_{n,\p})}
\newcommand{\fgnm}{\widehat{E}^-(K_{n,\p})}

\newcommand{\fgnje}{\widehat{E}^\ve(K_{n,\p_j})}

\newcommand{\fgnsubone}{\widehat{E}(K_{n-1,p})}

\newcommand{\LTpm}{L_p^\pm(E,\Kinf/K,X)}
\newcommand{\LTp}{L_p^+(E,\Kinf/K,X)}
\newcommand{\LTm}{L_p^-(E,\Kinf/K,X)}
\newcommand{\Lzne}{L_p^\ve(E,\Kinf/K,\zn - 1)}
\newcommand{\Lznp}{L_p^+(E,\Kinf/K,\zn - 1)}

\newcommand{\Lznpm}{L_p^\pm(E,\Kinf/K,\zn - 1)}
\newcommand{\Lpzero}{L_p^+(E,\Kinf/K,0)}
\newcommand{\Lmzero}{L_p^-(E,\Kinf/K,0)}

\begin{document}

\title{Iwasawa Theory of Elliptic Curves at Supersingular Primes over
$\Z_p$-extensions of Number Fields}

\author{Adrian Iovita and Robert Pollack}

\maketitle

\begin{abstract}
In this paper, we make a study of the Iwasawa theory of an elliptic curve at a supersingular prime $p$ along an arbitrary ${\bf Z}_p$-extension of a number field $K$ in the case when $p$ splits completely in $K$.  Generalizing work of Kobayashi \cite{Kobayashi02} and Perrin-Riou \cite{PR01}, we define restricted Selmer groups and $\lambda^\pm$, $\mu^\pm$-invariants; we then derive asymptotic formulas describing the growth of the Selmer group in terms of these invariants.  To be able to work with non-cyclotomic ${\bf Z}_p$-extensions, a new local result is proven that gives a complete description of the formal group of an elliptic curve at a supersingular prime along any ramified ${\bf Z}_p$-extension of ${\bf Q}_p$.
\end{abstract}


\section{Introduction}
\label{section:intro}

Over the last few years, much light has been shed on the subject of
Iwasawa theory of elliptic curves at supersingular primes.
In \cite{Kurihara01} and \cite{PR01}, asymptotic formulas for the size of
$\Sha(E/\Qn)[p^\infty]$ have been established where $\Qn$ runs through
the cyclotomic $\Zp$-extension of $\Q$.   In \cite{Kobayashi02} and
\cite{Pollack02}, a theory of algebraic and analytic $p$-adic
$L$-functions is formed that closely parallels the case of ordinary
reduction.  
The methods of all of the above papers depend heavily upon varying
the fields considered in the cyclotomic direction.  
The purpose of this paper is to extend some of these results to a
more general collection of $\Zp$-extensions.

The essential difference in Iwasawa theory between the ordinary  and
the supersingular case is that, in the later case, 
the Galois theory of Selmer groups is badly behaved.  Namely, if
$\Kinf/K$ is a $\Zp$-extension with layers $K_n$ and $E/K$ is an
elliptic curve supersingular at some prime over $p$, 
then the Selmer group of $E$ over $K_n$ is much smaller than
the $\Gal(\Kinf/K_n)$-invariants of the Selmer group of $E$ over
$\Kinf$.
(In the case of ordinary reduction, these two groups are nearly the same
by Mazur's control theorem.)
The reason descent fails in the supersingular case boils down to the fact
that the trace map on $\hE$ (the formal group of $E/\Qp$) is not
surjective along a ramified $\Zp$-extension.
Following \cite{PR90}, we make a careful 
study of how the trace map affects the Galois theory and we propose an
analogous ``control theorem'' 
that takes into account the formal group of
$E$ (see Theorem \ref{thm:control}).  In the end, this setup allows
one to convert the local 
information of $\hE$ into global information about the Selmer group of $E$.
These considerations are carried out in section \ref{section:control}.

In \cite{Kobayashi02}, a complete description of the Galois module
structure of $\hE(k_n)$ is given in terms of generators and relations
where $k_n$ runs through the local cyclotomic $\Zp$-extension of $\Qp$.  
The new local result of this paper is a generalization of the above
result to any ramified $\Zp$-extension of $\Qp$.
Namely, if $L_\infty/\Qp$ is a
ramified $\Zp$-extension with layers $L_n$ then we produce points $d_n \in
\hE(L_n)$ such that $\Trn(d_n) = -d_{n-2}$ for $n \geq 2$ where $\Trn :
\hE(L_n) \maps \hE(L_{n-1})$ is the trace map.  Furthermore, $d_n$ and
$d_{n-1}$ generate $\hE(L_n)$ over $\Zp[\Gal(L_n/\Qp)]$ (see Theorem
\ref{thm:formalgroups}).   From this result, we can completely
describe the kernel and cokernel of the trace map.
This local analysis is done in section \ref{section:formalgroup}.
Note that the above analysis of $\hE$ not only gives generators and
relations, but the generators satisfy a compatibility as the level
varies.  It is precisely this compatibility that allows Iwasawa
theory in the supersingular case to retain the flavor of the ordinary
case.  

To be able to apply these local results, we are obliged to work with
number fields $K$ for which $p$ splits completely (since the local
result assumes that we are working over $\Qp$).  For such $K$
and $p$, we analyze arbitrary $\Zp$-extensions of $K$.
Following \cite{PR90}, we produce algebraic $p$-adic
$L$-functions and then using the ideas of \cite{Kobayashi02} and
\cite{Pollack02} 
we form plus/minus $L$-functions that actually lie in the Iwasawa
algebra 
(assuming $a_p=0$).  Attached to these $L$-functions, we can 
associate plus/minus $\mu$ and $\lambda$-invariants.

In section \ref{section:kurihara}, we analyze the case where these
$L$-functions are units (i.e. when all the $\mu$ and
$\lambda$-invariants are zero).  In terms of $E$, this is the case
when $E(K)/pE(K) = 0$, $\Sha(E/K)[p]=0$ and $p \nmid \Tam(E/K)$; here 
$\Tam(E/K)$ represents the Tamagawa factor of $E$ over $K$.  (Note
that by the Birch and Swinnerton-Dyer conjecture, these hypotheses are
equivalent to $\frac{L(E/K,1)}{\Omega_{E/K}}$ being a $p$-adic unit.)
Under these strict global hypotheses, we prove that $E(K_n)$ and
$\Sha(E/K_n)[p^\infty]$ are finite for all $n$.  Furthermore, we
describe precisely the Galois structure of $\Sha(E/K_n)[p^\infty]$ and
in particular, produce precise formulas for its size.  

When $K=\Q$ and $\frac{L(E/\Q,1)}{\Omega_{E/\Q}}$ is a $p$-adic unit,
using Kato's Euler system we can verify our algebraic hypotheses and
we recover the main result of \cite{Kurihara01} (see
Corollary \ref{cor:kurihara}).  When $K$ is an 
imaginary quadratic extension of $\Q$ where $p$ splits and
$\frac{L(E/K,1)}{\Omega_{E/K}}$ is a $p$-adic unit, we can again
verify our algebraic hypotheses via Kato's result and produce exact
descriptions of the the size and structure of $\Sha(E/K_n)[p^\infty]$
(see Corollary \ref{cor:ac}).

In section \ref{section:plfunction}, we give two different constructions
of these plus/minus algebraic $p$-adic $L$-functions.
Namely, we follow \cite{PR90} and use the points $\{d_n\}$
to produce $p$-adic power series.  
Alternatively, we use the methods of \cite{Kobayashi02} to produce
restricted Selmer groups (which behave more like
Selmer groups at ordinary primes).  
These two approaches are related in that the characteristic power series
of the restricted Selmer groups agree with the power series constructed
(see Proposition \ref{prop:equiv}).

Finally, in section \ref{section:growth}, we study the arithmetic of $E$
along the extension $\Kinf/K$.  When the coranks of the Selmer groups grow without
bound along this extension, the algebraic $p$-adic $L$-functions vanish and the
restricted Selmer groups are not cotorsion (over the Iwasawa algebra).  
In this case, the coranks of
these restricted Selmer groups control the rate of growth of the coranks
of the Selmer groups at each finite level (see Proposition \ref{prop:rankgrowth}).
On the other hand, when these coranks remain bounded, we prove that
these $L$-functions are non-zero and the restricted Selmer groups are
indeed cotorsion.  In this case, we produce
asymptotic formulas for the growth of these Selmer group in terms of
the Iwasawa invariants of the plus/minus $L$-functions
as in \cite{PR01} (see Theorem \ref{thm:formulas}). 

~\\
\noindent
{\it Acknowledgments}: We are grateful to Ralph Greenberg for many 
interesting discussions on subjects pertaining to this paper.  We thank
Shin-ichi Kobayashi for conversations relating to \cite{Kobayashi02}.  We also
thank Nigel Byott and Cornelius Greither for very useful email exchanges.
Both authors were partially supported by NSF grants.

\section{Preliminaries}

Let $E/\Q$ be an elliptic curve and $p$ an odd supersingular prime 
for $E$.  Let $K/\Q$ be a finite extension
and $\Kinf/K$ a $\Zp$-extension with layers $K_n$.  Denote by 
$\Lambda$ the Iwasawa algebra $\Zp[[\Gal(\Kinf/K)]]$ 
and let $\Gamma_n = \Gal(\Kinf/K_n)$.
We will impose the following hypothesis on the splitting type of the 
prime $p$ in $\Kinf$.
\begin{itemize}
\item {\bf Hypothesis S} (for splitting type): The prime $p$ splits completely 
in $K$ 
into $d=[K:\Q]$ distinct primes,
say $\p_1,\dots,\p_d$.  Also, each $\p_i$ is totally ramified 
in $\Kinf$.
\end{itemize}

By abuse of notation, we will denote the unique prime over $\p_i$ in
either $K_n$ or $\Kinf$ simply by $\p_i$.  

\begin{lemma}
Hypothesis (S) implies that $E(K_{\p_i})[p]=0$ and $E(\Kinf)[p]=0$.
\end{lemma}

\begin{proof}
We have an exact sequence
\begin{eqnarray}
\label{eqn:reduction}
0 \maps E_1(\Qp) \maps E(\Qp) \maps \widetilde{E}(\Fp) \maps 0
\end{eqnarray}
where $\widetilde{E}$ denotes the reduction of $E$ mod $p$ and where $E_1$ is 
defined by the above sequence.  Since $p$ is supersingular,
$\widetilde{E}(\Fp)$ has no $p$-torsion (see \cite[V. Theorem 3.1]{Silverman}).  Furthermore, since
$E_1(\Qp) \cong \widehat{E}(\Qp)$, we have $E_1(\Qp)$ has no $p$-torsion
(see \cite[VII. Proposition 2.2 and IV. Theorem 6.1]{Silverman}).
Hence, $E(\Qp)[p]=0$.  Now since $p$ splits completely in 
$K$, we have that $K_{\p_i} \cong \Qp$ and $E(K_{\p_i})[p]=0$.  

For the second part, if $E(\Kinf)[p] \neq 0$ then $E(K)[p] =
E(\Kinf)[p]^{\Gamma} \neq 0$ since $\Gamma =\Gal(\Kinf/K)$ is pro-$p$.
However, $E(K)[p] \subseteq E(K_{\p_i})[p] = 0$.
\end{proof}

\subsection{Selmer groups}

For $L$ an algebraic extension of $\Q$ and $v$ a prime of $L$, define
$$
\H_E(L_v) = \frac{H^1(L_v,E[p^\infty])}{E(L_v) \otimes \Qp / \Zp} \text{~~and~~}
\P_E(L) = \prod_v \H_E(L_v)
$$
where the product is taken over all primes of $L$.  (Here $L_v$ is a union of 
completions of finite extensions of $\Q$ in $L$.)
Then the Selmer group of $E[p^\infty]$ is defined as
$$
\Sel(E[p^\infty]/L) = \ker \left( H^1(L,E[p^\infty]) \maps \P_E(L) \right).
$$
The Selmer group of $T_pE$ (the Tate module of $E$) is defined
similarly (the cocycles should locally lie in $E(L_v) \otimes \Zp$)
and will be denoted by $\Sel(T_pE/L)$.  We then have
\begin{eqnarray*}
0 \maps E(L) \otimes \Qp / \Zp \maps \Sel(E[p^\infty]/L) \maps
\Sha(E/L)[p^\infty] \maps 0
\end{eqnarray*}
and
\begin{eqnarray}
\label{eqn:ST}
0 \maps E(L) \otimes \Zp \maps \Sel(T_pE/L) \maps
T_p(\Sha(E/L)) \maps 0
\end{eqnarray}
where $\Sha(E/L)$ denotes the Tate-Shafarevich group and 
where $T_p(\Sha(E/L)) = \displaystyle \varprojlim_n \Sha(E/L)[p^n]$
is its Tate module (which is zero if $\Sha(E/L)$ is finite).

We will use the following abbreviations: 
$$
S_n = \Sel(E[p^\infty]/K_n),~S = \Sel(E[p^\infty]/\Kinf), 
$$
$$
S_n(T) = \Sel(T_pE/K_n),~
X_n = S_n^\dual~\text{and}~X=S^\dual
$$
where $Y^\dual = \Hom(Y,\Qp/\Zp)$.

\subsection{Local duality}

\begin{thm}[Tate Local Duality]
Let $v$ be a finite place of $K$.  There exists a perfect pairing
$$
H^1(K_v,E[p^\infty]) \times H^1(K_v,T_pE) \maps \Qp / \Zp
$$
induced by cup-product.
Furthermore, under this pairing $E(K_v) \otimes \Qp / \Zp$ is the
exact annihilator of $E(K_v) \otimes \Zp$, inducing an isomorphism
\begin{eqnarray}
\label{eqn:TLD}
\H_E(K_v)^\dual \cong E(K_v) \otimes \Zp.
\end{eqnarray}
\end{thm}

\begin{proof}
See \cite[Theorem 2.1]{Tate}.
\end{proof}

We can use local duality to analyze the local factors $\H_E(K_v)$
appearing in the definition of the Selmer group.
\begin{lemma}
\label{lemma:TLD}~
\begin{enumerate}
\item If $v \nmid p$ then $\H_E(K_v)$ is finite.  
\item If $\p | p$ then $\H_E(K_\p)^\dual \cong \widehat{E}(K_\p)$ 
assuming hypothesis (S).
\end{enumerate}
\end{lemma}

\begin{proof}
We have that $E(K_v) \cong \Z_l^{[K_v:\Qp]} \times T$ where $v \mid l$
and $T$ is a finite group (see \cite[VII. Proposition 6.3]{Silverman}).  Hence by (\ref{eqn:TLD}),
$\H_E(K_v) \cong \left( T \otimes \Zp \right)^\dual$ if $l \neq p$ and is therefore finite.  
For $\p | p$,
$$
\H_E(K_\p)^\dual \cong E(K_\p) \otimes \Zp \cong E_1(K_\p) \otimes \Zp
\cong \widehat{E}(K_\p)
$$
by (\ref{eqn:reduction}) since $p$ is supersingular and
$K_\p \cong \Qp$.
\end{proof}

\subsection{Global duality}

Let $\Sigma$ be a finite set of primes of $L$ containing $p$, the
infinite primes and all primes of bad reduction for $E$ and let
$K_\Sigma$ be the maximal extension of $K$ that is unramified outside
of $\Sigma$.  We have two exact sequences 
\begin{eqnarray}
\label{eqn:first}
0 \maps \Sn \maps H_n
\stackrel{\gamma_n}{\maps} \bigoplus_{v \in \Sigma} \H_E(K_{n,v})
\end{eqnarray}
and
\begin{eqnarray}
\label{eqn:second}
0 \maps \SnTstrict \maps \SnT \maps \bigoplus_{v \in \Sigma} E(K_v)
\otimes \Zp
\end{eqnarray}
where $H_n = H^1(K_\Sigma/K_n,E[p^\infty])$ and
where $\SnTstrict$ is defined by the second sequence.
By Tate local duality, $E(K_v) \otimes \Zp$
is dual to $\H_E(K_{n,v})$.  Global duality asserts that these two
sequences splice into a five term exact sequence.

\begin{thm}[Global duality]
The sequence
\begin{eqnarray*}
0 \maps \Sn \maps H_n
\stackrel{\gamma_n}{\maps} \bigoplus_{v \in \Sigma} \H_E(K_{n,v}) \maps
\SnT^\dual \maps \SnTstrict^\dual \maps 0
\end{eqnarray*}
is exact where the first two maps come from (\ref{eqn:first}) and the
last two maps come from (\ref{eqn:second}) and Tate local duality.
\end{thm}

\begin{proof}
For a statement of global duality in this form see \cite[Section 1.7]{Rubin}).
\end{proof}

\section{A control theorem in the supersingular case}
\label{section:control}

When $p$ is an ordinary
prime for $E$, Mazur proved that the natural map of
restriction between $\Sn$ and $\Sinfn$ has finite kernel and cokernel of size
bounded independent of $n$ (see \cite{Mazur72}).  
A theorem of this form, that compares $\Sn$ to $S^{\Gamma_n}$ is often
called a control theorem.
A key ingredient needed for this result is that
the trace map on the formal group of $E$ is surjective along a
ramified $\Zp$-extension. 

In the supersingular case, the trace fails to be surjective (see
\cite{Hazewinkel77}).  In fact, the $\Zp$-corank of 
$\coker\left( \Sn \maps \Sinfn \right)$ grows without bound.
In this section, we will produce an analogous control theorem that
describes this cokernel in terms of the formal group of $E$.
Throughout this section, we will be assuming (S).

Let $s_n : \Sn \maps \Sinfn$ and $r_{n,v} : \H_E(K_{n,v}) \maps
\H_E(K_{\infty,v'})$ denote the natural restriction maps with $v'$ some
prime of $\Kinf$ over $v$.
The following theorem can be thought of
as a control theorem in the supersingular case.  

\begin{thm}
\label{thm:control}
We have a four term exact sequence
$$
0 \maps \frac{\SnT}{\SnTstrict} \maps \fgn \times B_n \maps \Xinfn \stackrel{x_n}{\maps}
\Xn \maps 0 
$$
where $x_n = s_n^\dual$,
$\fgn = \oplus_{j=1}^d \widehat{E}(K_{n,\p_j})$ and
$B_n$ is a finite group whose size is bounded by 
the $p$-part of $\Tam(E/K_n)$.
\end{thm}

To prove this theorem, we will need to control the kernel and cokernel
of $s_n$.  We follow the methods of \cite{Greenberg:PC} and \cite{Greenberg:IEC}
and direct the reader to these articles for more details.

\begin{prop}
\label{prop:coker}
We have
\begin{enumerate}
\item $\ker(s_n) = 0$
\item $\coker(s_n) \cong \im(\gamma_n) \cap \left( \oplus_{v \in
  \Sigma} \ker(r_{n,v})\right)$
\end{enumerate}
where $\gamma_n$ is defined in (\ref{eqn:first}).
\end{prop}

\begin{proof}
This proposition follows from applying the snake lemma to the diagram
defining $\Sn$ and $S$.
See \cite[Chapter 4]{Greenberg:PC} especially Lemma 4.2 and 4.3 for details.
\end{proof}

The following proposition describes $\ker(r_{n,v})$.  The case of
primes dividing $p$ behaves quite differently from primes not over $p$.

\begin{prop}
We have
\label{prop:localker}
\begin{enumerate}
\item \label{item:1} For $v \nmid p$, $\ker(r_{n,v})$ is finite.  If
  $v$ splits completely in $\Kinf$ then $\ker(r_{n,v})=0$; otherwise
  it has size equal to $\Tam(E/K_{n,v})$ up to a $p$-adic unit.  
\item \label{item:2} For $\p | p$, $\ker(r_{n,\p}) = \H_E(K_{n,\p})$.
\end{enumerate}
\end{prop}

\begin{proof}
For part (\ref{item:1}), see the comments after Lemma 3.3 in
\cite{Greenberg:IEC}.
For part (\ref{item:2}), we have 
$$
\H_E(K_{\infty,\p}) = \varinjlim_n \H_E(K_{n,\p}) = \left( \varprojlim_n
\widehat{E}(K_{n,\p}) \right)^\dual
$$
where the last inverse limit is taken with respect to the trace map.
However, there are no universal norms for $\widehat{E}$ along the
ramified $\Zp$-extension $K_{\infty,\p}/K_\p$ since $p$ is
supersingular (see \cite{Hazewinkel77}).
Therefore, $\H_E(K_{\infty,\p}) = 0$ and $\ker(r_{n,\p}) = \H_E(K_{n,\p})$.
\end{proof}

\begin{remark}
The fact that $\ker(r_{n,\p})$ equals all of $\H_E(K_{n,\p})$
is the essential difference between the ordinary case and the
supersingular case and is the reason why the cokernel of $s_n$ grows
without bound.
\end{remark}

\begin{proof}[Proof of Theorem \ref{thm:control}]
To control $\coker(s_n)$
we will need to understand how $\im(\gamma_n)$ relates to 
$\oplus_{v \in \Sigma} \ker(r_{n,v})$.
To ease notation, let $\H_v = \oplus_{v \nmid p} \H_E(K_{n,v})$ and
$\H_p = \oplus_{\p | p} \H_E(K_{n,\p})$.
By global duality 
$$
\im(\gamma_n) = \ker \left( \H_v \times \H_p \maps
\SnT^\dual \right).
$$  
For $v \nmid p$, the image of 
$\H_E(K_{n,v})$ in $\SnT^\dual$ is zero (the former is a finite group
by Lemma \ref{lemma:TLD} and
the later is a free module by (S) and (\ref{eqn:ST}) ).  
Hence, we can write $\im(\gamma_n) = \H_v \times A$ with $A \subseteq
\H_p$ and applying global duality again yields
\begin{eqnarray}
\label{eqn:gd}
\frac{\H_p}{A} = \left( \frac{\SnT}{\SnTstrict} \right)^\dual.
\end{eqnarray}
By Proposition \ref{prop:localker}, $\oplus_{\p | p}
\ker(r_{n,\p}) = \H_p$ and hence 
$$
\im(\gamma_n) \cap \left( \oplus_{v \in \Sigma} \ker(r_{n,v}) \right) \cong
\left( \oplus_{v \nmid p} \ker(r_{n,v}) \right) \times A.
$$
Therefore, by Proposition \ref{prop:coker} and (\ref{eqn:gd}), we have
\begin{eqnarray}
\label{eqn:controldual}
0 \maps \coker(s_n) \maps \left( \oplus_{v \nmid p} \ker(r_{n,v}) \right)
\times \H_p \maps
\left( \frac{\SnT}{\SnTstrict} \right)^\dual \maps 0.
\end{eqnarray}
By Lemma \ref{lemma:TLD}, $\H_p^\dual \cong \fgn$
and by Proposition \ref{prop:localker}, $\# \ker(r_{n,v})$ is bounded
by the $p$-part of $\Tam(E/K_{n,v})$.
Therefore, dualizing (\ref{eqn:controldual}) yields the theorem.
\end{proof}

\section{Structure of some formal groups} 
\label{section:formalgroup}

\subsection{Lubin-Tate formal groups}
\label{sec:lubin_tate}

Let $p>2$ be a prime and  $\{L_n\}_{n\ge 0}$ with  $\displaystyle\Qp=L_0\subset L_1\subset L_2 ...\subset L_\infty=\cup_nL_n$
be a tower of fields such that $L_\infty$ is a totally ramified $\Zp$-extension of $\Q_p$. 
 Let $k_{n+1}:=L_n[\mu_p]$ for
$n\ge 0$ and
$k_\infty=\cup_n k_n=L_\infty[\mu_p]$. Here, if $M$ is a field by $M[\mu_p]$ we mean the 
extension of $M$ obtained by adjoining to $M$ the $p$-th roots of unity in some 
fixed algebraic closure of $M$. Then $k_\infty$ is a $\Zp^\times$-extension of $\Qp$ and 
the group of its universal norms is generated by a uniformizer
of $\Zp$, say $\pi$, such that
$\displaystyle \ord_p \left(\frac{\pi}{p}-1\right)>0$. Now
we would like to carefully choose
a Lubin-Tate formal group (by choosing a ``lift of Frobenius'' corresponding to $\pi$) whose
$\pi^n$-division points generate $k_n$ over $\Qp$. 

Namely, let us define 
$$
\displaystyle f(X):=\pi
X+\sum_{i=2}^p\frac{p(p-1)\cdots(p-i+1)}{i!}X^i\in \Zp[[X]].
$$
Then $f(X)$ is a lift of Frobenius corresponding to $\pi$, that is
$f(X)=\pi X \pmod{\text{deg~}2}$
and $f(X)=X^p \pmod{p}$ and moreover it satisfies the properties:
\begin{enumerate}
\item $f(X)=(X+1)^p-1 \pmod{p^2}$
\item the coefficient of $X^{p-1}$ is $p$.
\end{enumerate}
We call this  a {\bf good lift} of Frobenius. 

\begin{lemma}
\label{lemma:formal_mult}
For  $\pi$ as above and  for a  good lift of Frobenius $f(X)$ let us denote by $F_f(X,Y)$
the corresponding formal group law.  We have 
$F_f(X,Y)=X+Y+XY \pmod{p}$ and $[a]_f(X)=(X+1)^a-1 \pmod{p}$ for all $a\in \Zp$.
\end{lemma}

\begin{proof} 
Let us write $f(X)=(X+1)^p-1+p^2g(X)$ where $g(X)\in \Zp[[X]]$.
Then $F_f(X,Y)$ (respectively $[a]_f(X)$) is the unique power series with coefficients 
in $\Zp$ such that $F_f(X,Y)=X+Y \pmod{\text{deg~}2}$ and
$f(F_f(X,Y))=F_f(f(X),f(Y))$ (resp. such that $[a]_f(X)=aX
\pmod{\text{deg~}2}$ and $f([a]_f(X))=[a]_f(f(X))$). 

Writing the identity for $F_f$ we get: 
$$
(F_f(X,Y)+1)^p-1+p^2g(F_f(X,Y))=F_f((X+1)^p-1,(Y+1)^p-1)+p^2G(X,Y)
$$
for some  $G(X,Y) \in \Zp[[X,Y]]$.
Therefore $F_f(X,Y)$ satisfies the identity 
\begin{eqnarray}
\label{eqn:a1}
(F_f(X,Y)+1)^p-1=F_f((X+1)^p-1, (Y+1)^p-1)\pmod{p^2}.
\end{eqnarray}
If we write
\begin{eqnarray}
\label{eqn:a2}
 F_f(X,Y)=X+Y+\sum_{i,j\ge 1}a_{ij}X^iY^j
\end{eqnarray}
with $a_{ij}\in \Z/p^2\Z$, the coefficients $a_{ij}$ are obtained by
identifying the coefficients of the  
monomials of same degree in (\ref{eqn:a1}) and solving for $a_{ij}$.  
This process of solving for $a_{ij}$ requires division by a multiple
of $p$ (mod $p^2$) and therefore $a_{ij}$ is only uniquely
determined (mod $p$).
In other words, any power series as
in (\ref{eqn:a2}) above, satisfying (\ref{eqn:a1}),
is unique (mod $p$).  But the power series $X+Y+XY$ satisfies 
these conditions and
so we have $F_f(X,Y)=X+Y+XY \pmod{p}$. The proof for $[a]_f(X)$ is similar.
\end{proof}

Let us fix for the rest of this section $f(X)$ and $F_f(X,Y)$ as in Lemma 
\ref{lemma:formal_mult} and let 
$[i]_f(X)=\sum_{j=1}^\infty a_j(i)X^j$, for $i=1,2,...,p-1.$

\begin{cor} 
\label{cor:det}
We have that the determinant $\left(a_j(i)\right)_{1\le i\le p-1, 1\le j\le p-1}$ is  in $\Zp^\times$.
\end{cor}

\begin{proof} 
From Lemma \ref{lemma:formal_mult}, we see that 
$\left(a_j(i)\right)$ is a lower triangular matrix modulo $p$ with ones along
the diagonal.  Hence, $\det(a_i(j)) \in \Zp^\times$.
\end{proof}

Let us denote by $\cO_n$ the ring of integers in $k_n$ and by $M_n$
its maximal ideal. For every $n$ we have $k_n=\Qp[F_f[\pi^n]]$. 

\begin{cor}
\label{cor:result}
Let $\beta\in F_f[\pi^n]-F_f[\pi^{n-1}]$. Then  for every $1\le b\le
p-1$  we can find a linear combination with coefficients in $\Zp$ of
$[1]_f(\beta), [2]_f(\beta),...,[p-1]_f(\beta)$ which has the form  
$\beta^b+\beta^p V$ with $V\in \cO _n$. 


\end{cor}

\begin{proof} Apply Corollary \ref{cor:det}.
\end{proof}

For every $s\in \Z^{\ge 1}$ we denote by
$G_s(f)$ the $\Zp$-submodule of $M_s$ generated by $F_f[\pi^s]$.
The main result of this section is the following proposition.

\begin{prop}~
\label{prop:generation} 
\begin{enumerate}
\item \label{gen:parta} We have $G_n(f)=M_n$ for all $n\ge 1$.
\item \label{gen:partb} Each $\beta\in F_f[\pi^n]-F_f[\pi^{n-1}]$ generates $M_n/M_{n-1}$ as a 
$\Zp[\Gal(k_n/\Q_p)]$-module. 
\end{enumerate}
\end{prop}

\begin{proof} As $F_f[\pi^{n-1}]\subset M_{n-1}$ and as $\Gal(k_n/\Qp)(\beta)=F_f[\pi^n]-F_f[
\pi^{n-1}]$ part (\ref{gen:parta}) implies part (\ref{gen:partb}). 

To prove part (\ref{gen:parta}), let us first remark that $p\in G_s(f)$ for all
$s$, because Tr$_{k_s/k_{s-1}}(\beta)=-p$ for $\beta\in F_f[\pi^s]-F_f[\pi^{s-1}]$.
Hence in order to show that $G_n(f)=M_n$, we need to show that $G_n(f)$ contains 
elements of valuation $\displaystyle \frac{b}{p^n-p^{n-1}}$ for all $1\le b\le p^n-p^{n-1}-1$.
So far from Corollary \ref{cor:result} we know that $G_n(f)$ contains elements 
of valuation $\displaystyle \frac{b}{p^n-p^{n-1}}$ for all $0\le b\le p-1$.


\bigskip\noindent
Let us formulate the statement $P(s)$ depending on $s\ge 1$:

\bigskip\noindent 
$P(s)$: {\it For all $n\ge s$, all $\beta\in F_f[\pi^n]-F_f[\pi^{n-1}]$ and all 
$b$ with $0\le b \le p^{s}-p^{s-1}-1$,
there is a linear combination with coefficients in $\Zp$ of $[1]_f(\beta),[2]_f(\beta),...,
[p^s-p^{s-1}-1]_f(\beta)$ which has the form
$u\beta^b+\beta^{b+1}V$ where $u\in \cO_n^\times$ and $V\in \cO_n$.}

\bigskip
Notice that if for some $s\ge 1$, $P(s)$ is true then $G_s(f)=M_s$,
i.e. part (\ref{gen:parta}) is true for that $s$.
We will prove that the statement $P(s)$ is true for all $s\ge 1$ by induction on $s$.
The statement is true for $s=1$ by Corollary \ref{cor:result}.
 
Let us now suppose that $P(t)$ is true for all $1\le t< s$.
Let $n \geq s$, $\beta\in F_f[\pi^n]-F_f[\pi^{n-1}]$ and 
fix $b$ with $0\le b \le p^{s}-p^{s-1}-1$. 
Let us denote $\gamma:=[p]_f(\beta)$; then as $p=\pi u$ with $u\in
\Zp^\times$, $\gamma\in F_f[\pi^{n-1}]-F_f[\pi^{n-2}]$. 
Moreover, from Lemma \ref{lemma:formal_mult} it follows that we have
\begin{eqnarray}
\label{eqn:congGm}
[i+pj]_f(\beta)-[i]_f(\beta)-[j]_f(\gamma)=[i]_f(\beta)[j]_f(\gamma)+pU,
\end{eqnarray}
where $U\in \cO_n$, $0\le i\le p-1$ and $0\le j\le
p^{s-1}-p^{s-2}-1$. 

Now use the induction hypothesis on $\gamma$ to get that
a linear combination with coefficients in $\Zp$ of
$[1]_f(\gamma),[2]_f(\gamma),...,[p^{s-1}-p^{s-2}-1]_f(\gamma)$
has the form
$$
(v\gamma^c+\gamma^{c+1}W) +pU_1
$$
with $v\in \cO_{n-1}^\times$, $W\in \cO_{n-1}$ and $U_1\in \cO_n$. 
Then for each $i$ between $0$ and $p-1$, by (\ref{eqn:congGm}), there
is a linear combination with coefficients in $\Zp$ of the elements
$[1]_f(\beta),[2]_f(\beta),...,[p^{s}-p^{s-1}-1]_f(\beta)$
having the form
$$
[i]_f(\beta)(v\gamma^c+\gamma^{c+1}W) +pU_2
$$
with $U_2 \in \cO_n$.
Now again use Corollary \ref{cor:result} to get a linear combination
over $\Zp$ of $[1]_f(\beta),[2]_f(\beta),...,[p^{s}-p^{s-1}-1]_f(\beta)$
with the form
$$
(\beta^a+\beta^pT)(v\gamma^c+\gamma^{c+1}W)+pU_3
$$
with $ T, U_3\in \cO_n$.
We then have that a linear combination with coefficients in $\Zp$ of the desired elements has the form
$$
v\beta^a\gamma^c+\beta^p\gamma^{c}V_1=u\beta^b+\beta^{p(c+1)}V_2
$$
where $\displaystyle u=v\frac{\gamma^c}{\beta^{pc}}\in \cO_n^\times$ and $V_1,V_2\in \cO_n$.
This proves that $P(s)$ is true for all $s$ and hence $G_s(f)=M_s$ for all $s$.
\end{proof}

\subsection{Formal groups of elliptic curves with supersingular reduction}
\label{sec:ss}

Let $E/\Qp$ be an elliptic curve with supersingular reduction and suppose that $a_p=0$.  
Let us denote, as in the previous sections, by $\widehat{E}$ the formal group of $E$,
i.e. the formal scheme over $\Zp$  which is the formal completion of 
the N\'eron model of $E$ at the identity of its special fiber. Let $L_\infty / \Qp$ be
a ramified $\Zp$-extension with layers $L_n$. We denote by
$\Trn:\hE(L_n) \lra \hE(L_{n-1})$
the trace with respect to the group-law $\hE(X,Y)$. Then the following
is the main result of this section.

\begin{thm}~
\label{thm:formalgroups}
For $n \geq 0$ there exists $d_n \in \widehat{E}(L_n)$ such that
\begin{enumerate}
\item $\Tr^n_{n-1} d_n = - d_{n-2}$
\item $\Tr^1_0 d_1 = u \cdot d_0$ with $u \in \Zpx$.
\item For $n \geq 1$, $\widehat{E}(L_n)$ is generated by $d_n$ and $d_{n-1}$ as a
$\Zp[\Gal(L_n/\Qp)]$-module.  Also, $d_0$ generates $\widehat{E}(\Qp)$.
\end{enumerate}
\end{thm}

The proof of this theorem will fill the rest of this section. Let us consider the 
$\Zp^\times$-extension $k_\infty$ attached to $L_\infty$ as in the section \ref{sec:lubin_tate}
and denote by $\pi$ the generator of
the group of universal norms of the extension $k_\infty/\Qp$ which has positive valuation.
We will first construct a sequence of points $c_n\in \hE(k_n)$ which satisfy the same trace
relations. For this we will use Honda-theory as in section 8 of \cite{Kobayashi02} and we will
choose a particular representative of the isomorphism class of $\hE$ whose logarithm has a 
certain form. More precisely, let $f(X)$ be a ``good lift'' of Frobeius attached to
$\pi$ as in section \ref{sec:lubin_tate} and let
$$
\ell(X):=\sum_{k=0}^\infty (-1)^k\frac{f^{(2k)}(X)}{p^k}\in \Qp[[X]],
$$
where $f^{(0)}(X)=X$ and if $n\ge 1$ is an integer we set $f^{(n)}(X):=f(f^{(n-1)}(X))$.
By Honda theory, if we denote by $G(X,Y):=\ell^{-1}(\ell(X)+\ell(Y))$ then $\hE$ and $G$ 
are isomorphic formal groups over $\Zp$ and the logarithm of $G$ is
$\ell(X)$. For the rest of this section we will identify these two
formal groups and will write $\hE$ for $G$. 
We first have

\begin{lemma}
\label{lemma:torsion}
The formal group $\hE$ has no $p$-power torsion points in $k_n$ for all $n\ge 0$.  
\end{lemma}

\begin{proof} The proof is the same, modulo the obvious adjustments,
  as the proof of Proposition 8.8 of \cite{Kobayashi02}. 
\end{proof}

\begin{cor}
\label{cor:log_inj}
The group homomorphism $\ell :\hE(k_n)\lra \widehat{\bG}_a(k_n)$ is injective.
\end{cor}

\begin{proof}
This corollary follows immediately from Lemma \ref{lemma:torsion} since the kernel of the
logarithm of a formal group is composed precisely of the elements of finite order.
\end{proof}

Let $F_f(X,Y)$ be the Lubin-Tate formal group over $\Zp$ attached to the lift of Frobenius
$f(X)$ as in section \ref{sec:lubin_tate} and let us choose a $\pi$-sequence $\{e_n\}_{n\ge 0}$
in $k_\infty$, i.e. $e_n\in F_f[\pi^n]-F_f[\pi^{n-1}]$ such that $f(e_n)=e_{n-1}$ for all
$n\ge 1$. Let $\epsilon\in p\Zp$ be such that $\displaystyle \ell(\epsilon)=\frac{p}{p+1}$
and define $c_n\in \hE(k_n)$ to be $c_n=e_n[+]_{\hE}\epsilon$ for all $n\ge 0$. 
The following lemma computes the traces of the $c_n$.

\begin{lemma}
\label{lemma:c_n}
For $n\ge 1$, we have $\Trn(c_n)=-c_{n-2}$ where here $\Trn$ is the trace from $\hE(k_n)$ 
to $\hE(k_{n-1})$.  For $n=1$, $\Tr^1_0(c_1) = u \cdot c_0$ with $u \in \Zpx$.
\end{lemma}

\begin{proof} Everything is set up so that the proof follows formally
  the same steps as the proof 
of Lemma 8.10 in \cite{Kobayashi02}. Namely, as $\ell$ is injective on $\hE(k_\infty)$,
it is enough to show that the relation holds after applying $\ell$ to both sides of the equality.
For $n\ge 2$ we have
\begin{align*}
\ell(\Trn(c_n))&=\mbox{Tr}_{k_n/k_{n-1}}\left(\frac{p}{p+1}+\sum_{k=0}^\infty
 (-1)^k \frac{e_{n-2k}}{p^k} \right) \\
&=\frac{p^2}{p+1}-p+p\sum_{k=1}^\infty (-1)^k\frac{e_{n-2k}}{p^k}=-\ell(c_{n-2})
\end{align*}
where $e_k=0$ for $k$ negative.  The calculation is
similar for $n=1$. 
\end{proof}

\begin{prop}
\label{prop:technical}
We have $\ell(M_n)\subset M_n+k_{n-1}$ and $\ell$ induces an isomorphism
$$
\hE(k_n)/\hE(k_{n-1})\cong \ell(M_n)/\ell(M_{n-1})\cong M_n/M_{n-1}.
$$
In particular  $c_n$ generates 
$\hE(k_n)/\hE(k_{n-1})$ as a $\Zp[\Gal(k_n/\Qp)]$-module.
\end{prop}

\begin{proof} The proof follows the steps of the proof of Proposition 8.12
of \cite{Kobayashi02}. The main new ingredient is Proposition \ref{prop:generation}.
\end{proof}

\begin{cor}
\label{cor:gen_full}
For $n\geq1$, $c_n$ and $c_{n-1}$ generates $\hE(k_n)$ as a
$\Zp[\Gal(k_n/\Q_p)]$-module.  
\end{cor}

\begin{proof} This follows easily from Proposition \ref{prop:technical} and the trace 
relations satisfied by the $c_n$ (see Lemma \ref{lemma:c_n}).
\end{proof}

\begin{proof}[Proof of Theorem \ref{thm:formalgroups}]
Let $d_n:=\Tr_{k_{n+1}/L_n}(c_{n+1})\in \hE(L_n)$. Then it is easy to see that 
$\Trn(d_n)=-d_{n-2}$ for $n\ge 2$; so we only have to show that $d_n$
and $d_{n-1}$ generate $\hE(L_n)$ as a $\Zp[\Gal(L_n/\Q_p)]$-module
for $n\ge 1$. 
Let us denote by $\Delta$ the torsion subgroup of $\Zp^\times=\Gal(k_\infty/\Qp)$.
Then we have an isomorphism of $\Zp[\Gal(k_{n+1}/\Qp)]$-modules:
$$
\hE(k_{n+1})=\oplus_{\chi\in \widehat{\Delta}}\hE(k_{n+1})^\chi
$$
where $\widehat{\Delta}$ is the group of characters of $\Delta$ and $\hE(k_{n+1})^\chi$
is the maximal subgroup of $\hE(k_{n+1})$ on which $\delta\in \Delta$ acts by multiplication by
$\chi(\delta)$. The isomorphism is given by $x\mapsto (x^\chi)_{\chi\in 
\widehat{\Delta}}$ where $\displaystyle x^\chi:=\frac{1}{p-1}\sum_{\delta\in \Delta}\chi(\delta)
x^\delta$. Since $c_{n+1}$ and $c_n$ generate $\hE(k_{n+1})$ as a $\Zp[\Gal(k_{n+1}/\Qp)]$-module,
for every $\chi\in\widehat{\Delta}$, $c_{n+1}^\chi$ and $c_n^\chi$ generate $\hE(k_{n+1})^\chi$
over the same ring. In particular for $\chi$ equal to the trivial
character, we have $\hE(k_{n+1})^\chi= 
\hE(L_n)$, $(p-1)(c_{n+1}^\chi)=d_n$, $(p-1)(c_n^\chi)=d_{n-1}$ and the conclusion 
follows.
\end{proof}

The following proposition describes the relations that the $d_n$
satisfy.  We first introduce some notation that will be used
throughout the remainder of the paper.
Let $\displaystyle \Phi_n(X):=\sum_{i=0}^{p-1}X^{ip^{n-1}}$ be the
$n$-th cyclotomic polynomial, $\xi_n = \Phi_n(1+X)$ and $\wn(X):=(X+1)^{p^n}-1$.
Also set,
$$
\twnp:=\prod_{\begin{array}{c} 1\le m\le n \\ m \text{~even}\end{array}} \Phi_m(1+X),~~~~
\twnm:=\prod_{\begin{array}{c} 1\le m\le n \\ m \text{~odd}\end{array}} \Phi_m(1+X),~
$$
$\Xwnp = X \cdot \twnp$ and $\Xwnm = X \cdot \twnm$.  Note that $\wn = X
\cdot \twnp \cdot \twnm$.  Finally, set $\Ln = \Lambda / \wn \Lambda$.

\begin{prop}
\label{prop:rel}
There is an exact sequence
$$
0 \maps \hE(\Qp) \maps d_n \Ln \oplus d_{n-1} \Lnsubone \maps
\hE(L_n) \maps 0
$$
where the first map is the diagonal embedding (note that $\hE(\Qp) \subseteq
d_k \Lambda_k$ for each $k$) and the second map is $(a,b)
\mapsto a-b$.  Furthermore, $d_n \Ln \cong \Lambda / \Xwne$ with $\ve
= (-1)^n$.
\end{prop}

\begin{proof}
The exact sequence comes from Proposition 8.13 of \cite{Kobayashi02}.
For the second part, we have that
$$
\Xwne d_n = \Xwnsubtwoe (\xi_n \cdot d_n) = \Xwnsubtwoe
\Trn(d_n) = - \Xwnsubtwoe d_{n-2} = \dots = \pm X d_0 = 0.
$$
Hence, there is a surjective map $\Ln/\Xwne \maps d_n \Ln$ obtained by
sending $1$ to $d_n$.  To see that this map
is injective, it is enough to note that $\Ln/\Xwne$ and $d_n \Ln$ are
free $\Zp$-modules of the same rank (which follows from 
the above exact sequence).
\end{proof}

\begin{cor}
\label{cor:trace}
For $n \geq 0$ and $\ve = (-1)^n$,
\begin{enumerate}
\item $\kTn \cong  \Xwnsubtwoe d_n \Ln$.
\item $\ckTn$ is a $p$-group with $p$-rank equal to $q_n$ where
$$
\displaystyle 
q_n = \begin{cases} p^{n-1} - p^{n-2} + \dots + p - 1 &  2 | n \\
  p^{n-1} - p^{n-2} + \dots + p^2 - p & 2 \nmid n \end{cases}. 
$$
\end{enumerate}
\end{cor}

\begin{proof}
By Proposition \ref{prop:rel}, we have
$$
\label{eqn:reltrace}
\begin{CD}
0 @>>> \hE(\Qp) @>>> d_n \Ln \oplus d_{n-1} \Lnsubone @>>> \hE(L_n)
@>>> 0 \\
@. @VV\times pV @VVV @VV\Trn V \\
0 @>>> \hE(\Qp) @>>> d_{n-1} \Lnsubone \oplus d_{n-2} \Lnsubtwo @>>>
\hE(L_{n-1}) @>>> 0
\end{CD}
$$
where the middle vertical map sends $(d_n,0)$ to $(0,-d_{n-2})$ and
$(0,d_{n-1})$ to $(p \cdot d_{n-1},0)$.  Then applying the snake lemma and
Proposition \ref{prop:rel} yields the result.
\end{proof}

\subsection{The plus/minus Perrin-Riou map}

We follow closely section 8 of \cite{Kobayashi02} except that we work
with a $\Zp$-extension instead of a $\Zp^\times$-extension. This
produces a certain shift in the numbering but 
the main arguments are formally the same. Let
$T$ be the $p$-adic Tate-module of $E$ considered as a
$\Gal(\overline{\Q}_p/\Qp)$-module. 
The Kummer map $\hE(L_n)\lra H^1(L_n, T)$ together with cup product and the
Weil pairing induces
$$
(\ ,\ )_n:\hE(L_n)\times H^1(L_n, T)\lra H^2(L_n, \Zp(1))\cong \Zp.
$$
Let $G_n:=\Gal(L_n/\Qp)\cong \Z/p^n\Z$ and for every $x\in \hE(L_n)$ let us define
the morphism $P_{x,n}:H^1(L_n, T)\lra \Zp[G_n]$ by $\displaystyle P_{x,n}(z)=\sum_{\sigma\in G_n}
(x^\sigma, z)_n\sigma$. Both $H^1(L_n, T)$ and $\Zp[G_n]$ are naturally $G_n$-modules and
$P_{x,n}$ is $G_n$-equivariant for all $x$ and $n$. Moreover, for every
$x\in \hE(L_n)$ and $n\ge 1$ the following diagram 
$$
\begin{CD}
H^1(L_n,T) @>P_{x,n}>> \Z_p[G_n]\\
@VVV @VVV \\
H^1(L_{n-1},T) @>P_{\Trn(x),n-1}>>\Zp[G_{n-1}]
\end{CD}
$$
is commutative.
Using the sequence of points $\{d_n\}_n$ we consider two subsequences: $d_n^+=d_n$
if $n$ is even and $d_n^+=d_{n-1}$ if $n$ is odd and similarly $d_n^-=d_{n-1}$ if $n$ is even 
and $d_n^-=d_n$ if $n$ is odd. We set $\displaystyle
P_n^\pm:=(-1)^{[\frac{n+1}{2}]}P_{d_n^\pm,n}$ and define
$$
\displaystyle\hE^+(L_n):=\{ P\in \hE(L_n) \ | \  
\Tr^n_m(P)\in \hE(L_{m-1}) \text{~for~all~} 1 \le m\le n, m
\text{~odd}\}; 
$$
$$
\hE^-(L_n):=\{ P\in \hE(L_n) \ | \ 
\Tr^n_m(P)\in \hE(L_{m-1}) \text{~for~all~} 1\le m\le n, m
\text{~even} \}.
$$

\begin{lemma}
\label{lemma:E_pm}
$d_n^\pm$ generates $\hE^\pm(L_n)$ as a $\Zp[G_n]$-module.
\end{lemma}

\begin{proof} The proof is the same as the proof of Proposition 8.13
of \cite{Kobayashi02}. 
\end{proof}

We define $H^1_\pm(L_n, T):=\left(\hE(L_n)^\pm\otimes \Qp/\Zp\right)^\perp\subset
H^1(L_n, T)$ where we think of  
$\hE(L_n)^\pm\otimes \Qp/\Zp$ as embedded in $H^1(L_n, V/T)$ by the
Kummer map with $V=T\otimes_{\Zp}\Qp$. The orthogonal complement is
taken with respect to the Tate pairing $\displaystyle\langle \ ,\
\rangle :H^1(L_n, T)\times H^1(L_n, V/T)\lra \Qp/\Zp$.


\begin{lemma}~
\label{lemma:P_n}
\begin{enumerate}
\item $\ker(P_n^\pm)=H^1_\pm(L_n, T)$.

\item The image of $P_n^\pm$ is contained in $\twnmp \Ln$.
\end{enumerate}
\end{lemma}

\begin{proof} The first part is clear from Lemma \ref{lemma:E_pm}. 
For the second part, we have that
$$
\Xwnpm P_n^\pm(z) = \Xwnpm \sum_{\sigma \in G_n} \left(\left(d_n^\pm\right)^\sigma,z\right)_n
\sigma
= \sum_{\sigma \in G_n} \left(\Xwnpm \left(d_n^\pm\right)^\sigma,z\right)_n \sigma
= 0
$$ 
by Proposition \ref{prop:rel}.  The lemma then follows because any
element of $\Ln$ that is killed by $\Xwnpm$ is divisble by $\twnmp$.
\end{proof}


Since $\wn=X \twnp \twnm$, we have an isomorphism
$$
\Ln^\pm:=\Zp[X]/\Xwnpm\cong \twnmp \Ln.
$$ 
We define $P_{\Lambda,n}^\pm$ to be the unique map which makes the
following diagram commute.
$$
\begin{CD}
H^1(L_n, T) @>P_{\Lambda,n}^\pm>> \Ln^\pm\\
@VVV @VVV \\
\frac{H^1(L_n,T)}{H^1_\pm(L_n, T)} @>P_n^\pm>> \Ln
\end{CD}
$$
Here the right vertical map is $\Ln^\pm \cong \twnmp \Ln \subseteq \Ln$.
The properties of the maps $P_{\Lambda,n}^\pm$ are gathered in the
following proposition.

\begin{prop}~
\label{prop:prop}
\begin{enumerate}
\item For $n\ge 1$,
\begin{eqnarray}
\begin{CD}
\label{diag}
H^1(L_{n+1}, T) @>P_{\Lambda,n+1}^\pm>> \Lambda_{n+1}^\pm\\
@V \cores_{n+1/n}VV @VVV \\
H^1(L_n, T) @>P_{\Lambda,n}^\pm>>\Ln^\pm
\end{CD}
\end{eqnarray}
commutes.  (Here the right vertical map is the natural projection.)
\item $P^\pm_{\Lambda,n}$ is surjective for all $n\ge 1$.

\item \label{prop:c} $P_{\Lambda,n}^\pm$ determines an  isomorphism $\displaystyle 
\frac{H^1(L_n, T)}{H^1_\pm(L_n, T)}\cong \Lambda_n^\pm$.
\end{enumerate}

\end{prop}

\begin{proof} See the proofs of Proposition 8.22, 8.24 and 8.25 of \cite{Kobayashi02}.
\end{proof}

Diagram (\ref{diag}) allows us to consider the projective limit 
(with respect to $n$) of the maps $P_{\Lambda,n}^\pm$ and we denote
this limit by 
$$
\displaystyle
P^\pm_\Lambda:\bH^1(T):=\varprojlim_n (H^1(L_n, T), \cores)\lra
\varprojlim_n \Ln^\pm \cong \Lambda.
$$ 
Also, let $\displaystyle
\bH_\pm^1(T):=\varprojlim_n (H^1_\pm(L_n, T), \cores)$ 
and we have:

\begin{prop}
\label{prop:P_full}
$P_{\Lambda}^\pm$ defines an isomorphism $\bH^1(T)/\bH^1_\pm(T)\cong
\Lambda$.  Furthermore, $\bH_\pm^1(T)$ is a free $\Lambda$-module of
rank 1.
\end{prop}

\begin{proof} For the first part, 
apply part (\ref{prop:c}) of Proposition \ref{prop:prop}.  
Then, from the first part, we know that $\bH_\pm^1(T)$ is a
direct summand of $\bH^1(T)$ which by \cite[Proposition 3.2.1]{PR94} is a free
$\Lambda$-module of rank 2.  Therefore, $\bH_\pm^1(T)$ is a projective
$\Lambda$-module and since $\Lambda$ is local, $\bH_\pm^1(T)$ is free of rank 1.
\end{proof}

Finally, we have the following description of the maps $P_n^\pm$ in terms of the dual exponential
map of the Galois module $T$.

\begin{prop}
\label{prop:dual_exp}
We have
$$
P_n^\pm(z)=\left(\sum_{\sigma\in G_n}\ell(d_n^\pm)^\sigma\sigma\right)\left(\sum_{\sigma\in G_n}\exp^*_{\omega_E}
(z^\sigma)\sigma^{-1}\right).
$$
\end{prop} 

\begin{proof}
See Proposition 8.26 of \cite{Kobayashi02}.
\end{proof}

\section{The ``most basic'' case in Iwasawa theory}
\label{section:kurihara}

\subsection{Algebraic results}

In this section, we will be working under the following 
restrictive global hypothesis.  Recall that $p$ is assumed
to be odd.
\begin{itemize}
\item {\bf Hypothesis G} (for global): 
\begin{enumerate}
\item $p \nmid \Tam(E/K)$
\item $\Sha(E/K)[p] = 0$
\item $E(K)/pE(K) = 0$.
\end{enumerate}
\end{itemize}

In the good (non-anomalous) ordinary case, this hypothesis implies
that both the $\mu$-invariant and $\lambda$-invariant of $E$ vanishes
along any $\Zp$-extension of $K$.  For this reason, we refer to the
situation in this section as the ``most basic'' case.
Throughout this section we will be assuming (S) and (G) and under these
hypotheses we will prove the following theorem.

\begin{thm}
\label{thm:kurihara}
Assuming (SG), $a_p = 0$ and $p$ odd, we have
\begin{enumerate}
\item \label{part:1}
$E(K_n)$ is finite
\item \label{part:2}
$\Sha(E/K_n)[p^\infty]^\dual \cong \left( \Lambda / (\twnp,\twnm)
  \right)^d$
\item \label{part:3}
$\ord_p \left(\#\Sha(E/K_n)[p^\infty]\right)= d \cdot
  \sum_{k=0}^n q_k$
\end{enumerate}
where $d = [K:\Q]$ and 
$q_k$ is defined in Corollary \ref{cor:trace}.
\end{thm}

\begin{remark}
The hypothesis $a_p=0$ is probably not necessary.  See \cite{Pollack03} for
a proof of this theorem for general $a_p$ (divisible by $p$) when $K=\Q$.
However, the condition that $p$ is odd is necessary (see \cite[Remark 1.2]{Pollack03}).
\end{remark}

We begin by computing the structure of $X$ as a $\Lambda$-module.  The
following well known result does not assume (SG).

\begin{prop}
\label{prop:Lrank}
When $p$ is supersingular for $E/\Q$,
$$
\rk_{\Lambda} X = \cork_{\Lambda} H^1(K_{\Sigma}/\Kinf,E[p^\infty]) \geq d.
$$
\end{prop}

\begin{proof}
The first equality follows from \cite[Corollary 5]{Schneider79}.  The
inequality follows 
from a global Euler characteristic calculation 
(see \cite[Proposition 3]{Greenberg89}) since  
\begin{eqnarray*}
\cork_{\Lambda} H^1(K_{\Sigma}/\Kinf,E[p^\infty]) - \cork_{\Lambda}
H^2(K_{\Sigma}/\Kinf,E[p^\infty]) = d.
\end{eqnarray*}
\end{proof}

\begin{prop}
\label{prop:Xstruct}
Assuming (SG), X is a free $\Lambda$-module of rank $d$.
\end{prop}

\begin{proof}
Considering Theorem \ref{thm:control} with $n=0$ yields
\begin{eqnarray}
\label{eqn:base}
\widehat{E}(K_p) \surj X_\Gamma 
\end{eqnarray}
since, by (G), $\Sel(E[p^\infty]/K) = 0$ and $p \nmid \Tam(E/K)$.  
Now, $\rk_{\Zp} X_\Gamma \geq d$ by Proposition \ref{prop:Lrank} and
hence (\ref{eqn:base}) is an isomorphism since 
$\widehat{E}(K_p) \cong \Zp^d$.
By Nakayama's lemma, we can lift 
(\ref{eqn:base}) to a map $\Lambda^d \surj X$.  Again, by Proposition
\ref{prop:Lrank}, $\rk_\Lambda X \geq d$ and hence this map is an
isomorphism. 
\end{proof}

\begin{remark}
\label{rmk:gen}
For the remainder of this section we will fix an isomorphism of $X$
with $\Lambda^d$.  Such an isomorphism (as constructed in 
Proposition \ref{prop:Xstruct}) depends in part upon an
identification of $\widehat{E}(K_p)$ with $\Zp^d$.  We will now
specify this identification.  
By (S), $K_{\p_j} \cong \Qp$ and hence Theorem \ref{thm:formalgroups}
applies to $\widehat{E}(K_{n,\p_j})$.  Set $d_{n,j} =
(0,\dots,d_n,\dots,0) \in \fgn = \oplus_{i=1}^d
\widehat{E}(K_{n,\p_i})$ where $d_n \in 
\widehat{E}(K_{n,\p_j})$.  Then $\{d_{0,j}\}_{j=1}^d$ generates
$\widehat{E}(K_p)$ and in what follows we will assume that
$\widehat{E}(K_p)$ is identified with $\Zp^d$ via these generators.  
\end{remark}

In particular, Theorem \ref{thm:control} yields 
\begin{eqnarray}
\label{eqn:controlkurihara}
\fgn \stackrel{R_n}{\maps} \Lnd \maps \Xn \maps 0
\end{eqnarray}
where $\Ln = \Lambda / \wn \Lambda$.  Furthermore, for $m \leq n$ we have
\begin{eqnarray}
\label{diag:trace}
\begin{CD}
\fgn  @>R_n>>   \Lnd  \\
@V\Tr^n_mVV      @VVV \\
\fgm  @>R_m>>   \Lmd
\end{CD}
\end{eqnarray}
where $\Tr^n_m$ is the trace map and the right vertical map is the
natural projection.  
We postpone checking the commutativity of this diagram until section 6
(see Proposition \ref{prop:comm}).

\begin{lemma}
\label{prop:trivzero}
$\twnminuse | R_n(d_{n,j})$ with $\ve = (-1)^n$.
\end{lemma}

\begin{proof}

By Corollary \ref{cor:trace}, $d_{n,j}$ is killed by 
$\Xwne$.  Since $R_n$ is a Galois equivariant map,
$R_n(d_{n,j})$ is also killed by
$\Xwne$ and is therefore divisible by $\twnminuse$. 
\end{proof}

By Lemma \ref{prop:trivzero}, write 
$$
R_n(d_{n,j}) = \twnminuse \cdot (u_{1j}, \dots u_{dj}) \in \Lnd
$$
where $\ve = (-1)^n$.

\begin{lemma}
\label{lemma:det}
$\det(u_{ij})$ is a unit in $\Ln$.
\end{lemma}

\begin{proof}
To prove this lemma it is enough to check that
$\det(u_{ij}(0))$ is a unit in $\Zp$.  In the case that $n$ is even,
we have by diagram (\ref{diag:trace})
$$
R_n(d_{n,j}) \equiv R_0(\Tr^n_0(d_{n,j})) \text{~~in~~} \Lambda_0^d
\cong (\Zp[X]/X)^d.
$$
By Theorem \ref{thm:formalgroups}, $\Tr^n_0(d_{n,j}) = \pm
p^{\frac{n}{2}} d_{0,j}$.  Also by Remark \ref{rmk:gen}, we have
normalized $R_0$ so that $R_0(d_{0,j}) = (0,\dots,1,\dots,0)$
where 1 is in the $j$-th coordinate.  
Therefore, $R_n(d_{n,j})$ evaluated at $0$ equals 
$(0,\dots,\pm p^{\frac{n}{2}},\dots,0)$.

On the other hand,
\begin{align*}
R_n(d_{n,j})(0) &= \twnm(0) \cdot (u_{1j}(0), \dots , u_{dj}(0)) \\
                &= p^\frac{n}{2} \cdot (u_{1j}(0), \dots , u_{dj}(0)).
\end{align*}
Therefore, 
\begin{eqnarray}
\label{eqn:basegen}
u_{ij}(0) = 
\begin{cases}
0 & i \neq j \\
\pm 1 & i = j
\end{cases}
\end{eqnarray}
and $\det(u_{ij}(0)) = \pm 1 \in \Zpx$.
The case of $n$ odd is proven similarly using the fact that
$\Tr^1_0(d_{1,j}) = u \cdot d_{0,j}$ with $u \in \Zpx$.
\end{proof}

Let $I_n = R_n(\fgn) \subseteq \Lnd$.
Then by Theorem \ref{thm:formalgroups}, $I_n$ is 
the ideal of $\Lnd$ generated by $R_n(d_{n,j})$ and
$R_n(d_{n-1,j})$ for $j=1,\dots,d$.  

\begin{prop}
\label{prop:model}
$\Lnd / I_n \cong \left( \Lambda / (\twnp,\twnm) \right)^d$.
\end{prop}

\begin{proof}
Let $\twnje = (0, \dots , \twne,\dots,0)$ where $\twne$ lies in the $j$-th
coordinate and let $J_n$ be the ideal generated by $\twnjp$ and $\twnjm$
for $j = 1, \dots, d$.  To prove the proposition, it suffices to
show that $I_n = J_n$.  By Proposition \ref{prop:trivzero}, 
$I_n \subseteq J_n$. 
Conversely, from (\ref{eqn:basegen}) in the proof of Lemma
\ref{lemma:det}, we have that  
$(I_n)_\Gamma \cong (J_n)_\Gamma$.  Therefore, by Nakayama's lemma we
can conclude $I_n = J_n$.
\end{proof} 

\begin{proof}[Proof of Theorem \ref{thm:kurihara}]
From (\ref{eqn:controlkurihara}) and Proposition \ref{prop:model},
\begin{eqnarray}
\label{eqn:Xn}
\Xn \cong \left( \Lambda / (\twnp,\twnm) \right)^d.
\end{eqnarray}
An explicit computation (see \cite[Lemma 7.1]{Kurihara01}) shows that 
\begin{eqnarray}
\label{eqn:size}
\ord_p \# \left( \Lambda / (\twnp,\twnm) \right) = \sum_{k=0}^n q_k.
\end{eqnarray}
Therefore $\Sn$ is finite and in particular $E(K_n)$ is finite
proving part (\ref{part:1}).  Now since there is no presence of rank,
$\Sn \cong \Sha(E/K_n)[p^\infty]$; this together with (\ref{eqn:Xn})
yields part (\ref{part:2}). Finally, part (\ref{part:3}) follows from
(\ref{eqn:size}).
\end{proof}

\subsection{Analytic consequences}

We begin with a lemma that converts analytic hypotheses into algebraic
ones.  The following is a deep lemma that relies heavily upon Kato's
Euler system.

\begin{lemma}
\label{lemma:kato}
If $p$ is an odd supersingular prime for $E/\Q$ such that
\begin{enumerate}
\item $\ord_p \left( \frac{L(E/\Q,1)}{\Omega_{E/\Q}} \right) = 0$
\item $G_\Q \maps \Aut(E[p])$ is surjective
\end{enumerate}
then $\Sel(E[p^\infty]/\Q) = 0$ and $p \nmid \Tam(E/\Q)$.
\end{lemma}

\begin{proof}
We have that $L(E/\Q,1) \neq 0$ and hence from Kato's Euler system
\cite{Kato}, $E(\Q)$ and $\Sha(E/\Q)$ are both finite.  We must
show that $\Sha(E/\Q)[p^\infty]=0$ and $p \nmid \Tam(E/\Q)$.

The (analytic) $p$-adic $L$-function $L_p^{\text {an}}(E,T) \in
\Qpbar[[T]]$ interpolates special values of $L$-series and in
particular
$$
L_p^{\text {an}}(E,0) = \left(1 - \frac{1}{\alpha} \right)^2 
\frac{L(E/\Q,1)}{\Omega_{E/\Q}}
$$
where $\alpha$ is a root of $x^2 - a_p x + p$ (see \cite[Section 14]{MTT}).

In \cite{PR93}, Perrin-Riou constructed an algebraic $p$-adic
$L$-function $L_p^{\text {alg}}(E,T) \in \Qpbar[[T]]$ (defined up to a
unit in $\Lambda$) with the property that 
$$
L_p^{\text {alg}}(E,0) \sim \left(1 - \frac{1}{\alpha} \right)^2
\frac{\# \Sha(E/\Q) \cdot \Tam(E/\Q)}{\# E^{\text {tor}}(\Q)}
$$
when $\Sel(E[p^\infty]/\Q)$ is finite
(see also \cite[Th\'eor\`eme 2.2.1]{PR01}).

Kato proved a divisibility between these two $p$-adic $L$-functions
under the above assumption on the Galois representation.  Namely, we
have that
$$
L_p^{\text {alg}}(E,T) ~|~ L_p^{\text {an}}(E,T)
$$
in $\Zp[[T]]$ (see \cite[Theorem 12.5]{Kato} and 
\cite[Th\'eor\`eme 3.1.3]{PR01}).  In particular, 
$$
\ord_p \left( \Sha(E/\Q) \cdot \Tam(E/\Q) \right) \leq 
\ord_p \left( \frac{L(E/\Q,1)}{\Omega_{E/\Q}} \right).
$$
(Note that $E[p](\Q) = 0$ since $p$ is supersingular.)
From the above inequality, the lemma follows immediately since
we are assuming that the right hand side is zero.
\end{proof}

The following corollary, originally proven by Kurihara, follows from
Theorem \ref{thm:kurihara} and Lemma \ref{lemma:kato}.

\begin{cor}
\label{cor:kurihara}
Let $K=\Q$ so that $\Kinf = \Qinf$ is the cyclotomic $\Zp$-extension.
Let $E/\Q$ be an elliptic curve and $p$ an odd prime of good reduction with
$a_p=0$.  Assume that 
\begin{enumerate}
\item $\ord_p \left( \frac{L(E/\Q,1)}{\Omega_{E/\Q}} \right) = 0$
\item $G_\Q \maps \Aut(E[p])$ is surjective.
\end{enumerate}
Then the conclusions of Theorem \ref{thm:kurihara} hold with $d=1$.
\end{cor}

\begin{proof}
First note that $\ord_p \left( \frac{L(E/\Q,1)}{\Omega_{E/\Q}} \right) = 0$
implies that $\Sel(E[p^\infty]/\Q) = 0$ and
that $p \nmid \Tam(E/\Q)$ by Lemma \ref{lemma:kato}.
Therefore, hypothesis (G) is satisfied.  Furthermore, (S) is
automatically satisfied when $K=\Q$ and the conclusions of Theorem
\ref{thm:kurihara} follow.
\end{proof}

\begin{cor}
\label{cor:ac}
Let $K$ be a quadratic extension of $\Q$ and $\Kinf$ any
$\Zp$-extension of $K$. 
Let $E/\Q$ be an elliptic curve with $p$ an odd prime of good reduction
satisfying (S) for $K$ and such that $a_p=0$.  Assume further that 
\begin{enumerate}
\item $\ord_p \left( \frac{L(E/K,1)}{\Omega_{E/K}} \right) = 0$
\item $G_K \maps \Aut(E[p])$ is surjective.
\end{enumerate}
Then the conclusions of Theorem \ref{thm:kurihara} hold with $d=2$.
\end{cor}

\begin{proof}
Let $E^D$ be the quadratic twist of $E$ corresponding to $K/\Q$.  
Then $L(E/K,s) = L(E/\Q,s) \cdot L(E^D/\Q,s)$.  In particular,
$\ord_p \left( \frac{L(E/K,1)}{\Omega_{E/K}} \right) = 0$ implies
that $\ord_p \left( \frac{L(E/\Q,1)}{\Omega_{E/\Q}} \right) = 0$
and  $\ord_p \left( \frac{L(E^D,1)}{\Omega_{E^D/\Q}} \right) = 0$
since both special values are $p$-integral (see \cite[Remark 6.5]{Pollack02}).
Since $G_K$ surjects onto $\Aut(E[p])$, we have that $G_\Q$ surjects
onto both $\Aut(E[p])$ and $\Aut(E^D[p])$.  Therefore, by Lemma \ref{lemma:kato},
we have that $\Sel(E[p^\infty]/\Q) = \Sel(E^D[p^\infty]/\Q) = 0$ 
and that $p$ does not divide $\Tam(E/\Q) \cdot \Tam(E^D/\Q)$.  
From this we can conclude that
$\Sel(E[p^\infty]/K) = 0$ and that $p$ does not divide $\Tam(E/K)$.
Therefore, hypothesis (G) is satisfied and the conclusions of Theorem
\ref{thm:kurihara} follow.
\end{proof}

\section{Algebraic $p$-adic $L$-functions}
\label{section:plfunction}

In this section, we construct algebraic $p$-adic $L$-functions in two different 
ways.  First, we work directly with the points $\{d_n\}$ and Theorem \ref{thm:control}
to produce two $p$-adic power series as in \cite{PR90}.  However, as
in section \ref{section:kurihara}, we first remove certain trivial zeroes to obtain
elements of the Iwasawa algebra.  Alternatively, we consider plus/minus Selmer groups as
in \cite{Kobayashi02} and define algebraic $p$-adic $L$-functions as the characteristic
power series of these $\Lambda$-modules.  Finally, we show that these two constructions
yield the same power series (up to a unit in $\Lambda$).  We continue to assume (S)
in order to make use of the local results of section \ref{section:formalgroup}.

\subsection{Construction of algebraic $p$-adic $L$-functions via $\{d_n\}$}
\label{subsec:algpl}

We begin by generalizing the constructions done in section
\ref{section:kurihara}.  Assuming (G), it was shown in Proposition
\ref{prop:Xstruct} that $\rk_\Lambda X = d$.  In general, this would
be true assuming a form of the weak Leopoldt conjecture.  
We introduce this conjecture as another hypothesis.  (See \cite{Greenberg02}
for a formulation of this conjecture and for cases when it is known to
be true.)

\begin{itemize}
\item {\bf Hypothesis W} (for Weak Leopoldt):
 $\cork_\Lambda H^2(K_\Sigma / \Kinf,E[p^\infty]) = 0$.
\end{itemize}

\begin{prop}
\label{prop:Xrank}
When $p$ is supersingular for $E/\Q$, we have that (W) is equivalent
to $\rk_\Lambda X = d$.
\end{prop}

\begin{proof}
This is clear from Proposition \ref{prop:Lrank} and its proof.
\end{proof}

If $Y$ is the $\Lambda$-torsion submodule of $X$, we have
\begin{eqnarray}
\label{eqn:torsion}
0 \maps Y \maps X \maps Z \maps 0
\end{eqnarray}
where $Z$ is torsion free.  By Proposition
\ref{prop:Xrank}, embedding $Z$ into its reflexive hull yields a
sequence
\begin{eqnarray}
\label{eqn:hull}
0 \maps Z \maps \Lambda^d \maps H \maps 0
\end{eqnarray}
with $H$ finite.
We can then define a map
$$
\fgn \maps \fgn \times B_n \maps \Xinfn \maps Z_{\Gamma_n} \maps \Lnd
$$
where the second map comes from Theorem \ref{thm:control}, the third map
comes from (\ref{eqn:torsion}) and the final map comes from
(\ref{eqn:hull}).  Denote by $Q_n$ the map from $\fgn$ to
$\Xinfn$ and by $R_n$ the map from $\fgn$ to $\Lnd$.
These maps satisfy an important compatibility property already
exploited in section \ref{section:kurihara}.  Before discussing this
property, we state a lemma on the functoriality of the snake lemma.

\begin{lemma}
\label{lemma:snake}
For $i=1,2$, let
$$
\begin{CD}
@. A_i @>>> B_i @>>> C_i @>>> 0 \\
@. @Va_iVV @Vb_iVV @Vc_iVV @. \\
0@>>> A'_i @>>> B'_i @>>> C'_i @.
\end{CD}
$$
be a commutative diagram and assume that there are maps $A_1 \maps
A_2$, $B_1 \maps B_2$, 
$C_1 \maps C_2$ and likewise for $A'_i$, $B'_i$ and $C'_i$ such that
all the respective squares commute.  Then
$$
\begin{CD}
\ker(c_1) @>\delta_1>> \coker(a_1)\\
@VVV @VVV \\
\ker(c_2) @>\delta_2>> \coker(a_2)
\end{CD}
$$
commutes where $\delta_i$ is the boundary map coming from the snake lemma.
\end{lemma}

\begin{proof}
This follows from a diagram chase.
\end{proof}

\begin{prop}
\label{prop:comm}
For $m \leq n$, we have that the following diagrams
\begin{eqnarray*}
\begin{CD}
\fgn  @>Q_n>> \Xinfn @. ~~~~~~~~~@.  \fgn  @>R_n>>   \Lnd  \\
@V\Tr^n_mVV      @VVV  @. @V\Tr^n_mVV      @VVV \\
\fgm  @>Q_m>> \Xinfm  @.@. \fgm  @>R_m>>   \Lmd
\end{CD}
\end{eqnarray*}
commute.
\end{prop}

\begin{proof}
Since we have a fixed map $X \maps \Lambda^d$ defined independent of
$n$, the following square
$$
\begin{CD}
\Xinfn @>>> \Lnd \\
@VVV @VVV \\
\Xinfm @>>> \Lmd
\end{CD}
$$
commutes and therefore, we only need to check the commutativity of the
left diagram in the proposition. 

We will use the notation of Theorem \ref{thm:control}.
Furthermore, let $\hE_n = \hE(K_{n,p})$ and $K_n = \oplus_v
\ker(r_{n,v})$.  Then, examining the definition of $Q_n^\dual$ 
piece-by-piece yields
\addtolength{\arraycolsep}{-2pt}
$$
\begin{array}{ccccccccccc}
S^{\Gamma_m} &\surj &\coker(s_m)& \stackrel{\delta_m}{\cong} &K_m \cap
\im(\gamma_m) &\subseteq &K_m &\surj &\H_{m,p} &\cong &\left(\hE_m\right)^\dual \\

\da & &\da & &\da & &\da & &\da & &\da \\

S^{\Gamma_n} &\surj &\coker(s_n)& \stackrel{\delta_n}{\cong} &K_n \cap
\im(\gamma_n) &\subseteq &K_n &\surj &\H_{n,p} &\cong &\left(\hE_n\right)^\dual 
\end{array}
$$
where the first horizontal map (for either the top or bottom row) is
the natural projection, the second is given by the snake lemma
(Proposition \ref{prop:coker}), the third is the natural inclusion,
the fourth is the natural projection (applying Proposition
\ref{prop:localker}) and the fifth is given by Tate local duality.
The first vertical map is the natural inclusion, the second is induced
by this inclusion, the third, fourth and fifth maps are induced by
restriction and the sixth map is given by the dual of the trace map.

We now check the commutativity of this diagram square-by-square.  The
first square commutes essentially by definition.  The second square
commutes by the functoriality of the snake lemma (Lemma
\ref{lemma:snake}).  The third and 
fourth squares commute because restriction commutes with these natural
inclusions and projections.  Finally, the commutativity of the last
square is an essential property of Tate local duality 
(see \cite[Proposition 4.2]{Mazur72}). 
Dualizing then yields the proposition.
\end{proof}

Since these maps are Galois equivariant, 
Proposition \ref{prop:trivzero} remains valid in this setting.
In particular, we can write
$$
R_n(d_{n,j}) = \twnminuse \cdot (-1)^{\left[ \frac{n+1}{2} \right]}
\cdot (u_{1j}^n, \dots, u_{dj}^n)
$$
with $u_{ij}^n \in \Len$.

\begin{lemma}
\label{lemma:compatible}
For $n>1$ and $\ve=(-1)^n$, $u_{ij}^n \equiv u_{ij}^{n-2} \pmod{\Xwnsubtwoe}$.
\end{lemma}

\begin{proof}
This lemma follows from Theorem \ref{thm:formalgroups} and Proposition
\ref{prop:comm}.
\end{proof}

From Lemma \ref{lemma:compatible}, we have that $(u_{ij}^n)$ forms a
compatible sequence inside of $\displaystyle \varprojlim_n \Lminusen$ for
$n$ running through 
positive integers of a fixed parity.  When $n$ is even, denote this
sequence by $u_{ij}^+$ and when $n$ is odd by $u_{ij}^-$.  Since
$\displaystyle \varprojlim_n \Lpmn \cong \Lambda$, we can consider $u_{ij}^+$ and
$u_{ij}^-$ as Iwasawa functions.
We are now prepared to define the plus/minus algebraic $p$-adic
$L$-functions. 

\begin{defn}
Let $Y$ be the $\Lambda$-torsion submodule of $X$ and let $t_Y =
\char_\Lambda(Y)$.  Then set
$$
\LTpm := \det(u_{ij}^\pm) \cdot t_Y
$$
which is well-defined up to a unit in $\Lambda$.
\end{defn}

\begin{remark}
Note that $\LTpm$ can be identically zero.  This vanishing occurs when
$\cork_{\Zp}(S_n)$ is unbounded (see Corollary \ref{cor:nonzero}).
Furthermore, these coranks can indeed be unbounded.
For example, consider the case where $K$ is a
quadratic imaginary field and $\Kinf$ is the anticylotomic extension.
The recent results of \cite{Cornut02} show that if there are Heegner
points present then indeed the corank of $\Sn$ will grow without bound.
\end{remark}

\subsection{Restricted Selmer groups}
\label{section:pmselmer}

As in \cite{Kobayashi02}, we define plus/minus
Selmer groups by putting harsher local conditions at each $\p_i$.

\begin{defn}
Set
$$
\Sel^\pm(E[p^\infty]/K_n) = \ker \left( \Sel(E[p^\infty]/K_n) \maps
\prod_{\p | p} \frac{E(K_{n,\p}) \otimes \Qp / \Zp}{\hE^\pm(K_{n,\p}) \otimes \Qp / \Zp} \right)
$$
and $\displaystyle \Sel^\pm(E[p^\infty]/\Kinf) = \varinjlim_n \Sel^\pm(E[p^\infty]/K_n)$.
\end{defn}

These plus/minus Selmer groups behave like Selmer groups 
at ordinary primes.  In particular, they satisfy a control
theorem in the spirit of Mazur's original control theorem.

\begin{thm}
\label{thm:controlpm}
The natural map
$$
\Sel^\pm(E[p^\infty]/K_n)^{\Xwnpm = 0} \maps \Sel^\pm(E[p^\infty]/\Kinf)^{\Xwnpm = 0} 
$$
is injective and has a finite cokernel bounded independent of $n$.
\end{thm}

\begin{proof}
The proof in \cite[Theorem 9.3]{Kobayashi02} translates verbatim over
to our situation.
\end{proof}

If $X^\pm(E/\Kinf) = \Sel^\pm(E[p^\infty]/\Kinf)^\dual$ then $X^\pm(E/\Kinf)$ need not be
a torsion $\Lambda$-module.  The ranks of these modules will be discussed
in section \ref{subsec:coranks}.

\subsection{Comparing $\Sel^\pm(E[p^\infty]/\Kinf)$ and $\LTpm$}

As in the ordinary case, when $X^\pm(E/\Kinf)$ is a torsion module, its
characteristic power series should be considered as an algebraic $p$-adic
$L$-function.  The following proposition (whose proof will fill the
remainder of the section) relates this point of view with that of section
\ref{subsec:algpl}.

\begin{prop}
\label{prop:equiv}
Assuming (W),
$$
\char_\Lambda X^\pm(E/\Kinf) = \LTpm \cdot v.
$$
with $v \in \Lambda^\times$.
\end{prop}

Before proving this proposition, we begin with a few lemmas.

\begin{lemma}  
\label{lemma:wnkill}
For $\p | p$ and $n \geq 0$,
$$
\twnminuse \cdot \HneT \subseteq \fgne
$$ 
with $\ve = (-1)^n$.
\end{lemma}

\begin{proof}
We check this for $n$ even; the case of $n$ odd is similar.
We have
$$
\frac{\HnpmT}{\widehat{E}(K_{n,\p})} \inj
\frac{H^1(K_{n,\p},T)}{\HnmpT} \cong \Lambda / \Xwnmp \Lambda
$$
where the second map is given by Proposition \ref{prop:prop}.  
The first map is injective since $z \in \HnpT \cap \HnmT$ is
orthogonal to $\hE(K_{n,\p}) = \fgnp + \fgnm$ and hence in
$\hE(K_{n,\p})$.  However, this
map is not surjective; its image is killed by $\twnm$ (rather than
just $\wnm$) which we now check.

Note that $\frac{H^1(K_{n,\p},T)}{\HnmpT} \cong \Lambda /
\Xwnmp \Lambda$ is free over $\Zp$ and hence
$\frac{\HnpmT}{\widehat{E}(K_{n,\p})}$,
being a submodule, is also free.  Then since $\hE(K_{n,\p})$ is free, we
can conclude that $\HnpmT$ is free.  
Now 
$$
\rk_{\Zp} \frac{H^1(K_{n,\p},T)}{\HnpT} =
\rk_{\Zp} \Lambda / \Xwnp \Lambda = \deg(\Xwnp) = p^n - p^{n-1} + \dots
+ p^2 - p + 1.
$$
Hence, since $\rk_{\Zp} \HnT = 2 \cdot p^n$,
$$
\rk_{\Zp} \HnpT = 2 \cdot p^n -  (p^n - p^{n-1} + \dots
+ p^2 - p + 1) = p^n + q_n 
$$
and therefore $\frac{\HnpT}{\widehat{E}(K_{n,\p})}$ has $\Zp$-rank equal to $q_n$.

Now, note that any submodule of $\Lambda / \Xwnm
\Lambda \cong Z_p[X] / \Xwnm(X)$ of rank $q_n$ is of the form $p^rX
\Zp[X] / \Xwnm(X)$ and hence is annihilated by $\twnm$.  This proves
that
$$
\twnm \cdot \HnpT \subseteq \hE(K_{n,\p}).
$$
Since $\wnp \cdot \left( \twnm \cdot \HnpT \right) = \wn \cdot \HnpT = 0$, we further
have that $\twnm \cdot \HnpT \subseteq \hE^+(K_{n,\p})$ by the
definition of $\hE^+$; this  completes the proof.
\end{proof}

Repeating the arguments of Theorem \ref{thm:control} for the
plus/minus Selmer groups yields 
\begin{eqnarray}
\begin{CD}
B_n \times \bigoplus_{\p | p} \HnpmT  @>Q^\pm_n>>  
\Xinfn @>>> \Xnpm @>>> 0 \\
@AAA @AA=A @AAA \\
B_n \times \fgn  @>Q_n>> \Xinfn @>>> \Xn @>>> 0 \\ 
\end{CD}
\end{eqnarray}
where $\Xnpm = \Sel^\pm(E[p^\infty]/K_n)^\dual$.
Taking the projective limit of the top line of the above diagram yields
\begin{eqnarray}
\label{eqn:controlpm}
\bigoplus_{\p | p} \HinfpmT \stackrel{Q^\pm}{\maps} X \maps X^\pm(E/\Kinf) \maps 0.
\end{eqnarray}
Let $R^\pm$ be
the composition of $Q^\pm$ with the embedding of $X$ into $\Lambda^d$
from (\ref{eqn:torsion}) and (\ref{eqn:hull}) and define $R_n^\pm$
similarly.  
By Proposition \ref{prop:P_full}, $\HinfpmT$ is a free $\Lambda$-module of rank 1.

\begin{lemma}
\label{lemma:rel}
For each $j$, fix a generator $z_j$ of $\HinfjpmT$.  Then 
$$
R^\pm(z_j) = (u^\pm_{1j}, \dots, u^\pm_{dj}) \cdot v_j^\pm
$$
with $v_j^\pm$ a unit in $\Lambda$.  
\end{lemma}

\begin{proof}
Let $\ve=(-1)^n$.
We begin by recovering the sequence $\{d_{n,j}\}_n$ (constructed in
section \ref{section:formalgroup}) from the element $z_j$.  Let
$z^n_j$ be the image of $z_j$ in $\HnjpmT$ and let $d'_{n,j} =
(-1)^{\left[ \frac{n+1}{2} \right]} \cdot \twne \cdot z^n_j$.  

\bigskip
\noindent
{\bf Claim:} $d'_{n,j} = d_{n,j} \cdot v_j^\ve$ for $v^\ve_j$ a unit in $\Lambda$
(depending only on the parity of $n$).

\bigskip
\noindent
First note that by Lemma \ref{lemma:wnkill}, $d'_{n,j}$ is in fact an
element of $\fgnje$.  Furthermore, the $z^n_j$ are compatible under
corestriction by construction.  Therefore,
\begin{eqnarray}
\label{eqn:tracerel}
\Tr^n_{n-2}(d'_{n,j}) = \Tr^n_{n-2}(\twne z^n_j) = p \cdot
\twnsubtwoe z^{n-2}_j = -p \cdot d'_{n-2,j}.
\end{eqnarray}

Since $\fgne$ is cyclic, generated by $d_{n,j}$ (Lemma
\ref{lemma:E_pm}), we can write $d'_{n,j} = d_{n,j} \cdot v_{n,j}$ with $v_{n,j}
\in \Lambda / \Xwne \Lambda$.  Then (\ref{eqn:tracerel}) implies that
$(v_{n,j})_n$ forms a compatible sequence for $n$ of a fixed parity.
Call the limiting function in 
$\displaystyle \varprojlim_n \Lambda / \Xwne \Lambda \cong \Lambda$
by $v_j^+$ for $n$ even and by $v_j^-$ for $n$ odd.
To establish the claim it remains to show that $v_j^\pm$ is a unit.

By \cite[Proposition 9.2]{Kurihara01}, $\HinfjpmT$ surjects onto
$H^1_\pm(K_{\p_j},T) \cong \hE(K_{\p_j})$.  Therefore, $z^0_j$
(resp. $\Tr^1_0(z^1_j)$) generates $\hE(K_{\p_j})$.  In particular,
$d'_{0,j}$ (resp. $\Tr^1_0(d'_{1,j})$) differs from 
$d_{0,j}$ (resp. $\Tr^1_0(d_{1,j})$) by a unit in $\Zp$.  Hence,
$v^\pm_j(0) \in \Zpx$ and $v^\pm_j$ is a unit in $\Lambda$.

By the claim,
\begin{align*}
\twne \cdot R^\ve_n(z^n_j) &= R_n \left(
(-1)^{\left[\frac{n+1}{2}\right]}  d'_{n,j} \right) = 
R_n \left(
(-1)^{\left[\frac{n+1}{2}\right]}  d_{n,j} \right) \cdot v^\ve_j \\
&= \twne \cdot (u^n_{1j}, \dots , u^n_{dj}) \cdot v^\ve_j.
\end{align*}
Then cancelling $\twne$ and taking limits over $n$ of a fixed
parity yields the lemma.
\end{proof}

\begin{proof}[Proof of Proposition \ref{prop:equiv}]

If $\cork_{\Zp}(S_n)$ is unbounded we will see by 
Corollary \ref{cor:nonzero} and Corollary \ref{cor:torsion} that
our proposition holds with $(0)=(0)$.  So we may assume that
$\cork_{\Zp}(S_n)$ is bounded.
From (\ref{eqn:controlpm}), we have
\begin{align*}
\char_\Lambda X^\pm(E/\Kinf) &= \char_\Lambda \left( \frac{X}{Q^{\pm}\left(\bigoplus
  \HinfjpmT\right)} \right) \\ 
 &= \char_\Lambda \left( \frac{\Lambda^d}{\{R^\pm(z_j)\}_{j=1}^d}
  \right) \cdot \char_\Lambda Y 
\end{align*}
since $X/Y \subseteq \Lambda^d$ has finite index.  Then by Lemma
\ref{lemma:rel}, $\displaystyle\char_\Lambda \left(
\frac{\Lambda^d}{\{R^\pm(z_j)\}_{j=1}^d}   \right) =
\det(u^\pm_{ij}) \cdot v$ with $v = \prod_j v^\pm_j \in \Lambda^\times$.
Hence, $\char_\Lambda X^\pm(E/\Kinf) = \det(u^\pm_{ij}) \cdot \char_\Lambda Y
\cdot v = \LTpm \cdot v$ which completes the proof.
\end{proof}

\section{Growth of Selmer groups in $\Zp$-extensions}
\label{section:growth}

In this section, we explore the growth of $\cork_{\Zp}(S_n)$ as $n$ varies.
We describe this growth in terms of the $\Lambda$-ranks of $X^+(E/\Kinf)$ and $X^-(E/\Kinf)$.  
When $\cork_{\Zp}(S_n)$ is bounded, we compute the growth of $\Sha(E/K_n)[p^\infty]$
in terms of the $\mu$ and $\lambda$-invariants of $\LTpm$ as in \cite{PR01}.
Throughout this section, we will be assuming (S).

\subsection{Corank of Selmer groups}
\label{subsec:coranks}

Let $r^\pm = \rk_\Lambda X^\pm(E/\Kinf) = \cork_{\Lambda} \Sel^\pm(E[p^\infty]/\Kinf)$.

\begin{prop}
\label{prop:rankgrowth}
We have
$$
\cork_{\Zp}(S_n) = r^\ve \cdot q_n + r^{-\ve} \cdot q_{n-1} + \OO(1)
$$
where $\ve=(-1)^n$.  (Here, and in what follows, the $\OO(1)$ term depends upon
$E$ and upon $\Kinf/K$, but not upon $n$.)
\end{prop}

\begin{remark}
Note that if $r^+ = r^-$ then $\cork_{\Zp} S_n = r^\pm \cdot p^n + \OO(1)$
since $q_n + q_{n-1} = p^n - 1$.  In the ordinary case such growth formulas always have this 
form.  However, in the supersingular case, one should have situations where $r^+ \neq r^-$.
Namely, if $K$ is a quadratic imaginary extension of $\Q$, $\Kinf$ is the anti-cyclotomic
$\Zp$-extension and $E$ has CM by $K$, then conjecturally $r^\ve=1$ and $r^{-\ve}=0$ where
$\ve$ is minus the sign of the functional equation for $E$.  (See \cite[pg. 247]{Greenberg83})
\end{remark}

Before proving Proposition \ref{prop:rankgrowth}, we begin with a definition and some lemmas.

\begin{defn}
For $L$ a finite extension of $\Q$, let
$$
\Sel^0(E[p^\infty]/L) = \ker\left(\Sel(E[p^\infty]/L) \maps 
\prod_{\p | p} E(L_\p) \otimes \Qp / \Zp \right);
$$
$$
\Sel^1(E[p^\infty]/L) = \ker\left(\Sel(E[p^\infty]/L) \maps 
\prod_{\p | p} \frac{E(L_\p) \otimes \Qp / \Zp}{E(\Q) \otimes \Qp / \Zp} \right).
$$
\end{defn}
To ease notation, let $S_n^\pm = \Sel^\pm(E[p^\infty]/K_n)$, $S_n^0 = 
\Sel^0(E[p^\infty]/K_n)$, $S_n^1 = \Sel^1(E[p^\infty]/K_n)$ and $X^\pm =X^\pm(E/\Kinf)$.

\begin{lemma}
\label{lemma:parts}
We have
\begin{enumerate}
\item 
\label{part1}
For any $\Ln$-module $M$, the map $M^{\wne=0} + M^{\twnminuse=0} \maps M$
is surjective and has a finite kernel.

\item 
\label{part2}
$\left( \Sne\right)^{\twnminuse=0} \subseteq S_n^1$.

\item 
\label{part3}
The map $S_{n-1}^{-\ve} + S_n^1 \maps S_n^{-\ve}$ has finite cokernel and its
kernel is contained in $S_n^1$.

\item 
\label{part4}
The map $S_n^+ + S_n^- \maps S_n$ has finite cokernel and its kernel is 
contained in $S_n^1$.

\end{enumerate}

\end{lemma}

\begin{proof} Part (\ref{part1}) follows from the fact that $(\wne,\twnminuse)/(\wn)$
is finite.  To see part (\ref{part2}), note that if $\sigma \in \Sne$ is killed
by $\twnminuse$, then the restriction of $\sigma$ to any $\p$ over $p$ will lie in 
$$
\left( E^+(K_{n,\p}) \otimes \Qp/\Zp\right) 
\cap 
\left( E^-(K_{n,\p}) \otimes \Qp/\Zp \right) =
E(K_\p) \otimes \Qp/\Zp.
$$
Part (\ref{part3}) follows similarly.  Part (\ref{part4}) is 
\cite[Proposition 10.1]{Kobayashi02}.  
\end{proof}

\begin{lemma}
\label{lemma:sel1}
Assuming (W), $\cork_{\Zp} S_n^1$ is bounded independent of $n$.
\end{lemma}

\begin{proof}
Since $\coker\left(S_n^0 \maps S_n^1\right)$ has $\Zp$-rank bounded by $d$, it
suffices to check that $\cork_{\Zp} S_n^0$ is bounded.  
By Theorem \ref{thm:control}, we have that
$$
\rk_{\Zp} \SnT - \rk_{\Zp} \SnTstrict + \rk_{\Zp} X_{\Gamma_n} = \rk_{\Zp} \fgn + \rk_{\Zp} X_n.
$$
By (W), $\rk_\Lambda X = d$ and hence $\rk_{\Zp} X_{\Gamma_n} = d \cdot p^n + \OO(1)$.
Furthermore, we have that 
$$\
\rk_{\Zp} \SnT = \rk_{\Zp} X_n, ~~
\rk_{\Zp} \SnTstrict = \rk_{\Zp} S_n^0 \text{~~and~~}
\rk_{\Zp} \fgn = d \cdot p^n.  
$$
Hence,
$\cork_{\Zp} S_n^0$ is $\OO(1)$ (i.e. bounded).
\end{proof}

\begin{lemma}
\label{lemma:pm}
Assuming (W), for $\ve = (-1)^n$,
$$
\cork_{\Zp}\left( \Sne\right)^{\wne=0} +
\cork_{\Zp}\left( S_{n-1}^{-\ve} \right)^{\wnminuse=0} =
\cork_{\Zp} \Sn + \OO(1).
$$
\end{lemma}

\begin{proof}
By Lemma \ref{lemma:sel1} and part (\ref{part4}) of Lemma \ref{lemma:parts}, 
$$
\cork_{\Zp} \Sne + \cork_{\Zp} \Snme = \cork_{\Zp} \Sn + \OO(1).
$$ 
By part (\ref{part1}) and part (\ref{part2}) of Lemma \ref{lemma:parts}
and by Lemma \ref{lemma:sel1},
$$
\cork_{\Zp} \left( \Sne \right)^{\wne=0}
\cork_{\Zp} \left( \Snme \right)^{\wnminuse=0}
= \cork_{\Zp} \Sn + \OO(1).
$$
Finally, by Lemma \ref{lemma:sel1} and part (\ref{part3}) of Lemma \ref{lemma:parts}, 
$$
\cork_{\Zp} \left( \Sne \right)^{\wne=0}
\cork_{\Zp} \left( S_{n-1}^{-\ve} \right)^{\wnminuse=0}
= \cork_{\Zp} \Sn + \OO(1).
$$
\end{proof}

\begin{proof}[Proof of Proposition \ref{prop:rankgrowth}]
For any $m$, by Theorem \ref{thm:controlpm},
\begin{align*}
\cork_{\Zp}\left(\Sme\right)^{\wme=0} &= \rk_{\Zp}(X^\ve / \wme X^\ve) + \OO(1) \\
&= r^\ve \cdot \rk_{\Zp}\left(\Lambda / \wme \Lambda\right) + \OO(1) \\
&= r^\ve \cdot q_m + \OO(1).
\end{align*}
Taking $m=n$ and $n-1$, together with Lemma \ref{lemma:pm}, yields the proposition.
\end{proof}

When $\cork_{\Zp} S_n$ is bounded, we will see that $X^+$ and $X^-$ are
$\Lambda$-torsion.  We introduce this condition as another hypothesis.

\begin{itemize}
\item {\bf Hypothesis B} (for bounded): $\cork_{\Zp}(S_n)$ is bounded.
\end{itemize}

\begin{lemma}
(B) implies (W).
\end{lemma}

\begin{proof}
Let $r = \rk_{\Lambda} X$.  By Proposition \ref{prop:Lrank}
it suffices to check that $r=d$.  We have that $\rk_{\Zp} \fgn = d
p^n$ and by the theory of $\Lambda$-modules,
$\rk_{\Zp} \Xinfn$ grows like $r p^n$.  Then from Theorem
\ref{thm:control} and (B), we can conclude $r=d$ completing the proof.
\end{proof}

\begin{cor}
\label{cor:torsion}
Hypothesis (B) holds if and only if 
$X^+$ and $X^-$ are torsion $\Lambda$-module.
\end{cor}

\begin{proof}
If $X^+$ or $X^-$ are not torsion then by Theorem \ref{thm:controlpm},
(B) must fail.  Conversely, if (B) holds then (W) holds.  Hence, by
Proposition \ref{prop:rankgrowth}, $r^+ = r^- = 0$ and 
$X^\pm$ is torsion.
\end{proof}

\subsection{$p$-cyclotomic zeroes of $\LTpm$}

We now relate certain zeroes of 
the $p$-adic $L$-function $\LTpm$ to the corank of $S_n$.
Denote by $\zn$ a primitive $p^n$-th root of unity
and let $\xi_n = \Phi_n(1+X)$.

\begin{lemma}
\label{lemma:Ldfinite}
Let $u_j = (u_{1j},\dots,u_{dj}) \in \Ld$ for $j=1,\dots,d$.  Then we have
$\left( \Lambda / \xi_n \right)^d / (u_1,\dots,u_d)$ is finite if and
only if $\det(u_{ij}(\zn-1)) \neq 0$.  When this occurs 
$$
\ord_p \left( \# \frac{\left( \Lambda / \xi_n
  \right)^d}{(u_1,\dots,u_d)} \right) = \mu(f) \cdot (p^n-p^{n-1}) +
\lambda(f)  
$$
where $f = \det(u_{ij})$ and $n$ is sufficiently large.
\end{lemma}

\begin{proof}
Set $u_j(\zn-1) = (u_{1j}(\zn-1),\dots,u_{dj}(\zn-1))$ and then
$$
\frac{\left( \Lambda / \xi_n \right)^d}{(u_1,\dots,u_d)}
\cong \frac{\Zp[\mu_{p^n}]^d}{(u_1(\zn-1),\dots,u_d(\zn-1))}.
$$
By linear algebra, the left hand side is finite if and only if
$\det(u_{ij}(\zn-1)) \neq 0$.  Furthermore, when these groups are finite,
they have size $(p^n-p^{n-1}) \cdot \ord_p(\det(u_{ij}(\zn-1)))$
since $p$ is totally ramified in $\Zp[\mu_{p^n}]$.  Our result then
follows since for any non-zero $g \in \Lambda$, 
$$
\ord_p(g(\zn-1)) = \mu(g) + \frac{\lambda(g)}{(p^n-p^{n-1})}
$$
for $n$ large enough.
\end{proof}

\begin{prop}
\label{prop:infinite}
We have that 
\begin{enumerate}
\item $\Lpzero \neq 0$ and $\Lmzero \neq 0$ if and only if $S_0$ is finite.
\item For $n>1$, $\Lzne \neq 0$ for $\ve = (-1)^n$ if and only if
$\Sn/\Snsubone$ is finite.
\end{enumerate}
\end{prop}

\begin{proof}
We prove this for even $n>1$; the other cases follow similarly.
Consider the diagram
\begin{eqnarray}
\label{diag:levels}
\begin{CD}
\fgn  @>Q_n>>  \Xinfn @>>> \Xn @>>> 0  \\
@V\Tr^n_{n-1}VV      @VV\pi_nV  @VVV  \\
\fgnsubone  @>Q_{n-1}>> \Xinfnsubone @>>> \Xnsubone @>>> 0
\end{CD}
\end{eqnarray}
Then $\Sn/\Snsubone$ is finite if and only if $\ker(\pi_n) /
Q_n(\kTn)$ is finite by Corollary \ref{cor:trace}. We have that
$\ker(\pi_n) \cong \w_{n-1} X / \w_n X$ and from (\ref{eqn:torsion}) 
\begin{eqnarray}
\label{eqn:torsionwn}
0 \maps \frac{\w_{n-1} Y}{\w_n Y} \maps \frac{\w_{n-1} X}{\w_n X}
\maps \frac{\w_{n-1} Z}{\w_n Z} \maps 0.
\end{eqnarray}
The map $R_n$ restricted to $\kTn$ is given by the composite map
\begin{eqnarray}
\label{eqn:Pn}
\kTn \stackrel{Q_n}{\maps}  \frac{\w_{n-1} X}{\w_n X} \maps \frac{\w_{n-1} Z}{\w_n Z}
\maps \left( \frac{\w_{n-1} \Lambda}{\w_n \Lambda} \right)^d.
\end{eqnarray}
Now by Corollary \ref{cor:trace}, $\{ \Xwnsubtwop d_{n,j} \}_{j=1}^d$
generates $\kTn$ and we have that $R_n(\Xwnsubtwop d_{n,j}) = \wnsubone
\cdot (u^n_{1j},\dots,u^n_{dj})$.  Set $u^n_j =
(u^n_{1j},\dots,u^n_{dj}) \in \left( \Lambda / \xi_n \right)^d$.

First we consider the case where $t_Y(\zn-1) = 0$.  Then by
definition $\Lznp = 0$ and we need to check that $\Sn / \Snsubone$ is
infinite.  Since $t_Y(\zn-1) = 0$, we have that $ \w_{n-1} Y / \w_n Y$
is infinite.  Then by Proposition \ref{prop:Xrank} and
(\ref{eqn:torsionwn}) 
$$
\rk_{\Zp}\left( \wn X / \wnsubone X \right) > d \cdot (p^n -
p^{n-1}).
$$
But $\rk_{\Zp} \left( \kTn \right) = d \cdot (p^n -
p^{n-1})$ and hence $\Sn / \Snsubone$ is infinite from (\ref{diag:levels}).

So we may assume that $t_Y(\zn-1) \neq 0$.  Then 
$\Lznp \neq 0$ is equivalent to $\det(u^n_{ij}(\zn-1)) \neq 0$
which by Lemma \ref{lemma:Ldfinite} is equivalent to 
$\left( \Lambda / \xi_n \right)^d / (u^n_1,\dots,u^n_d)$ being finite.
Since the last two maps in (\ref{eqn:Pn}) have finite kernel and
cokernel, these last statements are equivalent to
$$
\frac{\wn X / \wnsubone X}{(Q_n(d_{n,1}),\dots,Q_n(d_{n,d}))}
= \frac{\ker(\pi_n)}{Q_n(\kTn)} \text{~~being~finite.}
$$
Then by (\ref{diag:levels}), this is equivalent to 
$\Sn / \Snsubone$ being finite completing the proof.
\end{proof}

\begin{cor}
\label{cor:nonzero}
$\LTp \neq 0 $ and $\LTm \neq 0$ if and only if $\cork_{\Zp}(\Sn)$ is bounded.
\end{cor}

\begin{proof}
The result follows from Proposition \ref{prop:infinite}
and the fact that a non-zero element of $\Lambda$ has finitely many
zeroes.   
\end{proof}

\subsection{Case of bounded rank}

Throughout this subsection, we will assume (B) and obtain formulas
describing the growth of $S_n$ along $\Kinf/K$.

\begin{defn}
Assuming (B) (so that $\LTpm$ is non-zero) define
$$
\lambda^\pm = \lambdapm = \lambda(\LTpm) 
$$
and
$$
\mu^\pm = \mupm = \mu(\LTpm).
$$
\end{defn}

We begin with a general lemma about the ``growth'' of
torsion $\Lambda$-modules.

\begin{lemma}
\label{lemma:stablewn}
If $Y$ is a torsion $\Lambda$-module then for $n$ large enough
$\wnsubone Y / \wn Y$ is finite of size $\mu(Y) \cdot
(p^n-p^{n-1}) + \lambda(Y) - \rk_{\Zp}(Y_{\Gamma_n})$. 
\end{lemma}

\begin{proof} 
By the structure theory of $\Lambda$-modules, we may assume that $Y$
is of the form $\Lambda / f^e$ with $f$ an irreducible polynomial.  
(Note that any finite groups that appear are killed by $\wn$ for $n$
large enough.)
If $\gcd(f,\wn)=1$ then 
$$
\frac{\wnsubone Y}{\wn Y} \cong \frac{Y}{\xi_n Y} \cong 
\frac{\Lambda}{(f^e,\xi_n)} \cong
\frac{\Zp[\mu_{p^n}]}{f^e(\zn-1)}.
$$
Now
\begin{align*}
\ord_p \left( \# \frac{\Zp[\mu_{p^n}]}{f^e(\zn-1)} \right) &= 
(p^n - p^{n-1}) \cdot \ord_p(f^e(\zn-1)) \\
&= (p^n - p^{n-1}) \cdot \mu(f^e) + \lambda(f^e)
\end{align*}
which implies the lemma since $\rk_{\Zp}(Y_{\Gamma_n})=0$.

If $f = \xi_k$ for some $k \leq n$ then
$$
\frac{\wnsubone Y}{\wn Y} \cong 
\frac{\Lambda}{(\xi_k^{e-1},\xi_n)} \cong
\frac{\Zp[\mu_{p^n}]}{\xi_k^{e-1}(\zn-1)}
$$
and the lemma follows since $\rk_{\Zp}(Y_{\Gamma_n}) = \deg(\xi_k)$.  
\end{proof}

The following duality theorem will be needed in what follows.  Let
$$
\SnTzero = \Sel^0(T_pE/K_n) = \ker \left( \SnT \maps \prod_{\p
  | p} E(K_{n,\p}) \otimes \Zp \right)
$$
so that $\SnTstrict \subseteq \SnTzero \subseteq \SnT$. 

\begin{thm}
\label{thm:duality}
Let $Y$ be the $\Lambda$-torsion submodule of $X$.  Then assuming (W),
$Y$ is pseudo-isomorphic to $\Sel_p^0(E[p^\infty]/\Kinf)^\dual$.
In particular, $\rk_{\Zp} Y_{\Gamma_n} = \rk_{\Zp} \SnTzero =
\rk_{\Zp} \SnTstrict$.
\end{thm}

\begin{proof}
The first statement is Corollary 2.5 in \cite{Wingberg89}.  For the
second statement, let $S_\infty^0 = \Sel_p^0(E[p^\infty]/\Kinf)$ and
$S_n^0 = \Sel_p^0(E[p^\infty]/K_n)$.  Then by \cite[Remark
  4.4]{Kurihara01}, $S^0_n$ and $\left(S^0_\infty \right)^{\Gamma_n}$
differ only by finite groups.  Therefore,
$$
\rk_{\Zp} Y_{\Gamma_n} = \cork_{\Zp} \left(S^0_\infty
\right)^{\Gamma_n}
= \cork_{\Zp} S^0_n = \rk_{\Zp} \SnTzero.
$$
Since, $E(K_{n,v}) \otimes \Zp$ is finite for $v \nmid p$, we have
that
$\rk_{\Zp} \SnTzero = \rk_{\Zp} \SnTstrict$ completing the proof.
\end{proof}

\begin{thm}
\label{thm:formulas}
Assuming (B), we have that
$$
\ord_p\left(\#(\Sn/\Snsubone)\right) = 
\begin{cases}
\mu^+ \cdot (p^n-p^{n-1}) + (\lambda^+ - s) \cdot n + d \cdot q_n & 2 | n \\
\mu^- \cdot (p^n-p^{n-1}) + (\lambda^- - s) \cdot n + d \cdot q_n & 2 \nmid n 
\end{cases}
$$
where $s$ is the stable value of $\cork_{\Zp} S_k$ and $n$ is
sufficiently large.
\end{thm}

\begin{proof}
Consider the diagram
\addtolength{\arraycolsep}{-2pt}
$$
\begin{array}{ccccccccccc}
0 & \maps &  \frac{\SnT}{\SnTstrict}  & \maps & B_n \times \fgn &
  \stackrel{Q_n}{\maps} &
\Xinfn  &\maps  & \Xn  & \maps  & 0 \\

 & & \downarrow & & \downarrow & & \downarrow & & \downarrow \\

0 & \maps &  \frac{S_{n-1}(T)}{S_{n-1,\Sigma}(T)}  & \maps & B_{n-1} \times \fgnsubone &
  \stackrel{Q_{n-1}}{\maps} &
\Xinfnsubone  &\maps  & \Xnsubone  & \maps  & 0
\end{array}
$$
defined by Theorem \ref{thm:control}.
For $n$ large enough, $\SnT,\SnTstrict$ and $B_n$ stabilize and the
vertical maps in the above diagram between these groups become
multiplication by $p$.

We will break the above diagram into two pieces; namely
\begin{eqnarray}
\label{diag:left}
\begin{CD}
0 @>>> \frac{\SnT}{\SnTstrict} @>>> B_n \times \fgn  @>>>  M_n @>>>0 \\
@.   @V\cdot pVV  @V\cdot p \times \Tr^n_{n-1}VV      @VVm_nV  @.  \\
0 @>>> \frac{S_{n-1}(T)}{S_{n-1,\Sigma}(T)} @>>> B_{n-1} \times \fgnsubone
@>>> M_{n-1} @>>> 0
\end{CD}
\end{eqnarray}
and
\begin{eqnarray}
\label{diag:right}
\begin{CD}
0 @>>> M_n @>Q_n>>  \Xinfn @>>> \Xn @>>> 0  \\
@.   @Vm_nVV  @VVV      @VV\pi_nV  @.  \\
0 @>>> M_{n-1} @>Q_{n-1}>> \Xinfnsubone @>>> \Xnsubone @>>> 0
\end{CD}
\end{eqnarray}
where $M_n$ is defined by the above diagrams.  If $s_0 = \rk_{\Zp}
\SnTstrict$ and $h = \rk_{\Fp}B_n/pB_n$ then applying the snake lemma to
(\ref{diag:left}) yields
\begin{align*}
0 \maps &(\Z/p\Z)^h \times \kTn \maps \ker(m_n) \maps (\Z/p\Z)^{s-s_0}
\maps \\ 
&(\Z/p\Z)^h \times \ckTn \maps \coker(m_n) \maps 0.
\end{align*}
By Corollary \ref{cor:trace}, $\ckTn \cong (\Z/p\Z)^{d q_n}$ and
therefore we have that
\begin{eqnarray}
\label{eqn:kerrn}
\ker(m_n) \cong \kTn \times (\Z/p\Z)^{h+a} 
\end{eqnarray}
and
\begin{eqnarray}
\label{eqn:cokerrn}
\coker(m_n) \cong (\Z/p\Z)^{d q_n+h-s+s_0+a} 
\end{eqnarray}
for some $a$ between $0$ and $s-s_0$.
Applying the snake lemma to (\ref{diag:right}) yields
\begin{eqnarray}
\label{eqn:snakeright}
0 \maps \ker(m_n) \maps \frac{\wnsubone X}{\wn X} \maps \left( \Sn /
\Snsubone \right)^\dual \maps \coker(m_n) \maps 0.
\end{eqnarray}

For $n$ large enough, $\Sn/\Snsubone$ and 
$\frac{\wnsubone Y}{\wn Y}$ are both finite and
\begin{align*}
\left[ \frac{\wnsubone X}{\wn X} : \kTn \right] 
&= \# \left( \frac{\wnsubone Y}{\wn Y} \right) \cdot 
\left[ \frac{\wnsubone Z}{\wn Z}: \kTn \right] \\
&= \# \left( \frac{\wnsubone Y}{\wn Y} \right) \cdot
\left[ \frac{\wnsubone \Lambda}{\wn \Lambda}: \kTn \right] \\
&= \# \left( \frac{\wnsubone Y}{\wn Y} \right) \cdot
\# \left( \frac{(\Lambda/\xi_n)^d}{(u^n_1,\dots,u^n_d)}  \right).
\end{align*} 
Again, for $n$ large enough,  $\Lznpm \neq 0$ 
and hence $\det(u^n_{ij}(\zn-1)) \neq 0$.  Also,
by Lemma \ref{lemma:Ldfinite},
\begin{eqnarray}
\label{eqn:Ldsize}
\ord_p \left( \# \frac{\left( \Lambda / \xi_n
  \right)^d}{(u^n_1,\dots,u^n_d)} \right) = 
(\mu^\ve - \mu_t) \cdot (p^n-p^{n-1}) + \lambda^\ve - \lambda_t
\end{eqnarray}
where $\lambda_t = \lambda(t_Y)$, $\mu_t = \mu(t_Y)$ and $\ve = (-1)^n$.
Then, by Lemma \ref{lemma:stablewn}, we have that
$$\ord_p \left( \# \frac{\wnsubone Y}{\wn Y} \right) =
\mu_t \cdot (p^n-p^{n-1}) + \lambda_t - \rk_{\Zp}(Y_{\Gamma_n})$$
for $n$ large enough.  Thus,
$$
\ord_p \left[ \frac{\wnsubone X}{\wn X} : \kTn \right] 
= \mu^\ve \cdot (p^n-p^{n-1}) + \lambda^\ve - \rk_{\Zp}(Y_{\Gamma_n}).
$$
Returning to  (\ref{eqn:snakeright}), we can compute
\begin{align*}
\ord_p \left( \# \Sn / \Snsubone \right) 
&= \ord_p \left[ \frac{\wnsubone X}{\wn X} : \ker(m_n) \right]+ \ord_p(\# \coker(m_n)) \\
&= -a-h + \ord_p \left[ \frac{\wnsubone X}{\wn X} : \kTn \right] +
\ord_p(\# \coker(m_n)) \\
&=  \mu^\ve \cdot (p^n-p^{n-1}) + \lambda^\ve - \rk_{\Zp}(Y_{\Gamma_n}) + 
dq_n  - s +s_0 
\end{align*}
Finally, from Theorem \ref{thm:duality}, we have that $s_0 = \rk_{\Zp}
(Y_{\Gamma_n})$ which completes the proof of the theorem.
\end{proof}

\end{document}